\documentstyle[12pt]{article}
\makeatletter
\renewcommand{\@cite}[2]{[{{\bf #1}\if@tempswa,#2\fi}]}
\renewcommand{\@biblabel}[1]{[#1]\hfill}
\makeatother

\newtheorem{prop}{Proposition}
\newtheorem{nt}{Remark}
\newtheorem{th}{Theorem}
\newtheorem{lemma}{Lemma}

\newfont{\sdbl}{msbm9}
\newfont{\dbl}{msbm10 at 12pt}
\newcommand{\eqdef}{\stackrel{\rm def}{=}}
\newcommand{\proof}{{\bf Proof.\ }}

\newcommand{\lto}{\longrightarrow}
\newcommand{\oo}{{\cal O}}

\newcommand{\f}{{\cal F}}
\newcommand{\g}{{\cal G}}
\newcommand{\h}{{\cal H}}
\newcommand{\U}{{\cal U}}
\newcommand{\A}{{\cal A}}
\newcommand{\V}{{\cal V}}

\newcommand{\Ker}{{\rm Ker}\:}
\newcommand{\Spec}{{\rm Spec}\:}
\newcommand{\Coker}{{\rm Coker}\:}
\newcommand{\Image}{{\rm Im}\:}
\newcommand{\Mor}{{\rm Mor}\:}
\newcommand{\Ob}{{\rm Ob}\:}

\newcommand{\lm}{\mathop{\rm lim}}
\newcommand{\dm}{\mathop{\rm dim}}
\newcommand{\il}{\lm\limits_{\longrightarrow}}
\newcommand{\pl}{\lm\limits_{\longleftarrow}}
\begin{document}
\author{D. V. Osipov}
\title{Krichever correspondence for algebraic varieties}
\date{}
\maketitle
\begin{abstract}
In the work is constructed
new acyclic resolutions of quasicoherent sheaves.
These resolutions is connected with
multidimensional local fields.
Then the obtained resolutions is applied
for a construction of generalization of
the Krichever map to algebraic varieties
of any dimension.

This map  gives in the canonical way
two $k$-subspaces
$B \subset k((z_1)) \ldots
((z_n))$ and  $W \subset k((z_1)) \ldots ((z_n))^{\oplus r}$
from arbitrary
algebraic $n$-dimensional
Cohen-Macaulay  projective integral scheme $X$  over a field $k$,
a  flag of closed integral subschemes
$X=Y_0 \supset
Y_1 \supset \ldots Y_n$
(such that $Y_i$
is an ample Cartier divisor on
 $Y_{i-1}$, and $Y_n$
is a smooth $k$-point on all $Y_i$),
formal local parameters of this flag in the point
 $Y_n$,
a rank $r$  vector  bundle  $\f$ on  $X$,
and a trivialization $\f$
in the formal neighbourhood of the point
$Y_n$,
where the $n$-dimensional local field $k((z_1)) \ldots ((z_n))$
is associated with the flag
$Y_0 \supset \ldots \supset Y_n$.
In addition,
the constructed map is injective,
i.~e., it is possible to reconstruct uniquely all
the original geometrical data.
Besides,
from the subspace $B$ is written explicitly a complex,
which calculates  cohomology of the sheaf
$\oo_X$ on $X$;
and from the  subspace $W$
is written explicitly a complex,
which calculates cohomology of
$\f$ on $X$.
\end{abstract}

\section{Introduction}
In 70's years I.~M. Krichever suggested a construction
how to attach to some algebraic-geometric data,
connected with algebraic curves and vector bundles on them,
an infinite-dimensional (Fredholm) subspace
in the space $k((z))$ of Laurent power series (~\cite{Krich}).
This construction was successfully used
in the theory of integrable systems,
in paricular, in the theory of KP and KdV equations
(\cite{Krich}, \cite{SW}, \cite{D}).

 There were  also found
applications of this construction
to the theory   of modules of algebraic curves (\cite{ADKP}, \cite{BS}).
Besides, this construction turned out
to be connected  with description of commutative
subrings in the rings of pseudo-differential operators (\cite{M}, \cite{D}).
Now this construction is called the Krichever correspondence
or the Krichever map (\cite{ADKP},
\cite{M}, \cite{BF}, \cite{Q}).
But in these works it is essentially
that algebraic-geometric data
are connected with $1$-dimensional varieties
and $1$-dimensional local field $k((z))$.

Recently, in works~\cite{Par}, \cite{P} it
 were pointed out  some connections
between the theory of the KP-equations
and $n$-dimensional local fields;
also it was suggested
a variant of the Krichever map
for algebraic-geometric data
which  is connected
with algebraic surfaces,
vector bundles on them
and $2$-dimensional local fields.

One of the typical examples
of multidimensional local field
is the field of Laurent iterated series
 $k((z_1)) \ldots ((z_n))$.

Such fields serve
for natural generalization of
local objects of $1$-dimensional varieties
to the case of multidimensional varieties.
Let us consider an $n$-dimensional algebraic scheme  $X$.
Let  $Y_0 \supset  \ldots
  \supset Y_n$ be a flag of closed subschemes on
 $X$ such that
 $Y_0 = X$, $Y_i$ is of codimension $1$ in $Y_{i-1}$, and
$Y_n=x$ is a closed point.
Then there exists a construction
(\cite{P2}, \cite{Bejl}, \cite{PL}),
attaching in the canonical way to such flag  some ring,
which is an $n$-dimensional local field
provided that $x$ is a smooth point on all $Y_i$.
Moreover,
if
$X$
is an algebraic variety over a field $k$,
$x$ is a  $k$-rational point,
and we fix  local parameters $z_1, z_2, \ldots, z_n
\in \hat{\cal O}_{x, X}$ such that $z_{n-i+1}=0$ is a
local equation of  variety $Y_i$
in the formal neighbourhood of the point  $x$
on the variety $Y_{i-1}$ ($1 \le i \le n$),
then the obtained $n$-dimensional local field
it is possible to identify with $k((z_1))((z_2)) \ldots ((z_n))$.

A concept  of multidimensional local field
 has  appeared
in the middle of 70's years,
and, originally, such fields were used
for the development of generalization
of class field theory
to the schemes of higher dimension.
Later there were also found
applications
of multidimensional local fields
to many problems of algebraic geometry,
where it makes sence to speak
about local components of geometric objects (see~\cite{PF}).

In this work, using multidimensional local fields,
we construct the Krichever cor\-res\-pon\-den\-ce
for varieties  of arbitrary dimension $n$:
that is some injective map from algebraic-geometric
data, connected with projective algebraic varieties,
full flags  of ample divisors and their local parameters
in the formal neighbourhood of the last point of the flag,
 vector bundles and their trivialisations in the formal
neighbourhood of the last point of the flag,
to some $k$-subspaces
of finite dimensional vector space over
the  $n$-dimensional local field $k((z_1)) \ldots ((z_n))$.

If $n=1$,
then our constructed map is
a variant of the Krichever map for curves.

If $n=2$,
then  our constructed map coincides
with the map constructed in~\cite{P}.

\medskip
The work  is organized as follows.

In \S2
we give various technical lemmas
about cohomology of coherent sheaves,
projective and injective limits,
which will be useful further in the work.

In \S3
we give a construction
of  family of functors,
which is connected with
quasicoherent sheaves
and a fixed flag of subvarieties,
and which can be interpreted
as a cohomology system of coeficients
on the standart symplex.

In \S4,
using the construction of \S3,
we construct complexes of sheaves of abelian groups,
which is acyclic resolutions of arbitrary
quasicoherent sheaves on schemes.

In \S5
we prove some theorems about
intersections among components of resolutions
constructed in \S4.
In some cases the whole resolution
can be reconstructed
from one $k$-subspace of finite dimensional
vector space over
 $k((z_1)) \ldots ((z_n))$.
Using this, we construct
the Krichever map in higher dimensions.

Note that  connected with
$n$-dimensional local fileds
resolutions of quasicoherent sheaves
on schemes
were in works~\cite{Bejl},
\cite{H2}.
But in contrast to these works,
our resolutions depend only on a single flag of subvarieties
and are not resolutions of adelic type.

Note also that,
in contrast to~\cite{P},
all the constructions and proofs in this work
are internal ones, i.~e., they are not  reduced
to  multidimensional  adelic complexes.

\bigskip
\medskip

During  all the work
we shall keep the following
notations and agreements.

For any finite set
$I$ let $\sharp I$ be
the number of elements
of the set $I$.

If  $X$  is a scheme, then\\
$Sh(X)$ is the category
of sheaves of abelian groups on $X$, \\
$CS(X)$   is
the category of coherent sheaves on
 $X$, \\
$QS(X)$     is
the catehory of quasicoherent sheaves on $X$,\\
$Ab$  is  the catehory of abelian groups.

If $f : Y \longrightarrow X$ is a morphism of two schemes,
then always  $f^*$ is
the pull-back   functor  in the category
of sheaves of abelian groups, $f_*$ is
the direct image functor
in the category of sheaves of abelian groups.

If $\U$ is an open covering of $X$,
$\f$ is a sheaf of abelian groups on $X$,
then $\check{H}^*(\U, \f)$~ are
the \v{C}ech cohomologies groups with respect
to the covering $\U$.

Let $Y \hookrightarrow X$ be a closed subscheme of a scheme $X$,
which is defined by the ideal sheaf $J$. Then by $(Y, \oo_X/J^k)$
denote the scheme
whose topological space coincides with the topological
space of the scheme
$Y$
and the structure sheaf is
 $\oo_X/J^k$. ($\oo_X$ is the structure sheaf of the scheme $X$.)

\bigskip
The author would like
to express the deep gratitude
to his scientific adviser
A.~N.~Parshin
for the constant attention to
the work.

\section{Technical lemmas}
\begin{lemma}   \label{lem1}
On a noetherian scheme $X$
any short exact sequence of
quasicoherent sheaves
is direct limit
of short exact sequences of coherent sheaves.
For $\phi : \f \longrightarrow \g \; \in \Mor(QS(X))$
there are $\phi_i : \f_i \longrightarrow \g_i \; \in \Mor(CS(X))$
with $\lm\limits_{\longrightarrow} \phi_i = \phi$,
$\f_i \subset \f$, $\g_i \subset \g$.
\end{lemma}
\proof See~\cite[ lemma 1.2.2]{H2} and \cite[ lemma 2.1.5]{H1}.

\medskip

\begin{lemma}     \label{lem2}
Let $X$ be a noetherian scheme.
Let $\psi : CS(X) \longrightarrow Sh(X)$ be
an exact additive functor.
Then $\psi$
commutes with direct limits.
\end{lemma}
\proof
(By analogy with
\cite[ lemma 1.2.3]{H2}
 or
\cite[ lemma 2.2.2]{H1}.)

First let us prove that
if we have a direct system of sheaves
$$
\left\{ \f_i : i \in I \; , \; \phi_{ij} : \f_i \to \f_j (i \le j)
        \right\}
$$
         with  $\lm\limits_{\longrightarrow} \f_i = 0 $,
then $\lm\limits_{\longrightarrow} \psi(\f_i) = 0 $.

For this one we prove that
for any open
$U \subset X$
$\lm\limits_{\longrightarrow} H^0(U, \psi(\f_i))=0 $.
Let $x \in \lm\limits_{\longrightarrow} H^0(U, \psi(\f_i)) $.
Let this
$x$
be represented by
 $x_i \in H^0(U, \psi(\f_i))$.
Then from   coherent property of the sheaf
$\f_i$
and  noetherian property of the scheme $X$
there is some $j \in I$ such that $\phi_{ij} = 0$.
Since $\psi$ is an additive functor,
we have that  $\psi (\phi_{ij}) = 0$.
Therefore
$$
\begin{array}{rcl}
^0(U, \psi(\f_i)) & \longrightarrow & H^0 (U, \psi(\phi_{ij}) (\psi(\f_i))) \\
x_i & \mapsto & 0
\end{array}
$$
Now consider the general case:
let $ \lm\limits_{\longrightarrow} \f_i = \f $,
$\phi_i : \f_i \to \f$ be the canonical morphisms.
Consider the following exact sequence of coherent sheaves:
$$
0 \longrightarrow \Ker \phi_i \lto \f_i \lto \f \lto \Coker \phi_i \lto 0
$$
The functor $\psi$ is an exact functor,
therefore we have the following exact sequence:
$$
0 \longrightarrow \psi(\Ker \phi_i) \lto \psi(\f_i) \lto
\psi(\f) \lto \psi(\Coker \phi_i) \lto 0
$$
From  $\il \Ker \phi_i =0$
and
$\il \Coker \phi_i =0$
it follows by arguments above that
$\il \psi ( \Ker \phi_i) =0$
and
$\il \psi (\Coker \phi_i) =0$.
Direct limit maps exact sequences to exact sequence.
Therefore $\il \psi (\f_i) = \psi(\f)$. Lemma~\ref{lem2} is proved.

\medskip

\begin{lemma}     \label{lem3}
Let $X$ be a noetherian scheme.
Then an exact additive functor
$\psi : CS(X) \lto Sh(X)$
can be uniquely extended
to a functor
$\psi' : QS(X) \lto Sh(X) $
which commutes with direct limits.
This new functor is exact as well.
\end{lemma}
\proof (By analogy with \cite[ lemma 1.2.4]{H2}.)

Let $\f \in \Ob(QS(X))$.
By lemma~\ref{lem1}
$\f = \il \f_i$,
where $\f_i \in \Ob(CS(X)) $.
Define
$$
\psi'(\f)= \il \psi(\f_i)  \mbox{.}
$$
We have
 $\psi' (\f) = \psi(\f) $
for $\f \in Ob (CS(X))$ by lemma~\ref{lem2}.
By lemma~\ref{lem1},
for any $\phi \in \Mor (QS(X))$
we have
$\phi = \il \phi_i$,
where $\phi_i  \in \Mor (CS(X))$.
Define
$$
\psi'(\phi) = \il \psi (\phi_i)      \mbox{.}
$$
By lemma~\ref{lem2}, we have that
 $\psi'(\phi) = \psi(\phi)$
for $\phi \in \Mor (CS(X))$.
It is clear that
this definition is the only one possible.
And by lemma~\ref{lem2},
it is well defined.
Lemma~\ref{lem3} is proved.

\medskip

\begin{lemma}  \label{lem5}
Let $X$ be a noetherian scheme,
$i : Y \hookrightarrow X$ be
a closed subscheme, which is defined by  the ideal sheaf
$J$ on $X$.
Let
$j : U \hookrightarrow Y$ be
an open subscheme of $Y$ such that
for any point
 $x \in X$
there exists an affine neighborhood
 $V \ni x$
such that $V \cap U$ is an affine subscheme.
Let the supports of sheaves $\f_i \in Sh(X)$  ($i=1, \ldots, 3$)
are in $Y$,
and the sheaf $\f_1$ is a quasicoherent sheaf
with respect to the subscheme
 $(Y, \oo_X / J^k)$
for some $k \in N$.
Then from  exactness of the sequence of sheaves
\begin{equation}  \label{eqn2}
0 \lto \f_1 \lto \f_2  \lto \f_3 \lto 0
\end{equation}
it follows  exactness of the following sequence
$$
0 \lto i_*j_*j^* \f_1 \lto    i_*j_*j^* \f_2
\lto i_*j_*j^* \f_3 \lto 0   \quad \mbox{.}
$$
\end{lemma}
\proof
First note that for any affine open
subscheme
$W \subset U$
and for any quasicoherent
sheaf
$\g$ on the scheme $(U, \oo_X/J^k \mid_U)$
we have
\begin{equation}  \label{eqn3}
H^1 (W, \g) = 0  \quad \mbox{.}
\end{equation}
In fact,
if the sheaf $\g$ is a quasicoherent sheaf
with respect to the subscheme $\;$
$U = (U, \oo_X/ J \mid_U)$,
then equality~(\ref{eqn3})
follows from  affineness of the scheme
$W$.
Now if
$\f$ is  a quasicoherent sheaf
with respect to the subscheme
$(U, \oo_X / J^k \mid_U)$,
$k \in {\bf N}$, $k \ge 1$,
then consider the following exact sequence:
\begin{equation}   \label{eqn4}
0 \lto J\f \lto \f \lto \f / J\f \lto 0   \quad \mbox{.}
\end{equation}
But the sheaves $J \f$
and $\f / J \f$
are quasicoherent sheaves
with respect to the subscheme
$ (U, \oo_X / J^{k-1} \mid_U)  $.
Therefore we can do induction,
from which it follows
that
$$
H^1(W, J\f) = 0  \qquad \qquad \mbox{and}
\qquad
\qquad H^1(W, \f/ J \f)=0  \; \mbox{.}
$$
Hence and from the long cohomological sequence
associated with sequence~(\ref{eqn4})
we obtain equality~(\ref{eqn3}).

Return to sequence~(\ref{eqn2}).
We have  exactness of the following sequence:
$$
0 \lto j^* \f_1 \lto j^* \f_2 \lto j^* \f_3 \lto 0 \quad \mbox{.}
$$
Appliing the functor
$j_*$,
we obtain
$$
0 \lto j_*j^* \f_1 \lto j_*j^* \f_2  \lto j_* j^* \f_3 \lto R^1j_* (j^* \f_1)
$$
Let us show that the sheaf
$$
R^1j_* (j^* \f_1) = 0 \quad \mbox{.}
$$
The sheaf
$j^* \f_1$
is a quasicoherent sheaf with respect
to the subscheme $(U, \oo_X / J^k \mid_U)$
for some $k \in {\bf{N}}$,
therefore the sheaf
$R^1j_* (j^* \f_1)$
is a quasicoherent
sheaf on the  scheme
$(Y, \oo_X/J^k)$
with respect to  the same
$k \in {\bf{N}}$.
Therefore it suffices to show that
 for affine open
$V$
from the lemma's conditions
\begin{equation} \label{eqn5}
 H^0 (V \cap Y, R^1 j_* (j^* \f_1)) = 0  \quad \mbox{.}
\end{equation}
But
$ H^0 (V \cap Y, R^1 j_* (j^*\f_1))
= H^1 (V \cap U, j^* \f_1)
$.
And equality~(\ref{eqn5})
follows from equality~(\ref{eqn3}).
Therefore we have  exactness of the following sequence:
$$
0 \lto j_* j^* \f_1 \lto j_* j^* \f_2  \lto j_* j^* \f_3  \lto 0  \quad \mbox{.}
$$
From
$i: Y \hookrightarrow X$ is a closed
imbedding it follows that
$i_*$ is an exact functor.
Therefore the following sequence is exact:
$$
0 \lto i_*j_*j^* \f_1 \lto    i_*j_*j^* \f_2
\lto i_*j_*j^* \f_3 \lto 0   \quad \mbox{.}
$$
Lemma~\ref{lem5} is proved.

\medskip
\bigskip
\bigskip
\bigskip

Let $X$ be a noetherian  scheme.
Suppose that we have an exact and additive functor
$\Phi \: : \: QS(X) \lto Sh(X)$.
Let $i: Y \hookrightarrow X$ be
a closed subscheme of the scheme
 $X$,
which is defined by the ideal  sheaf $J$.
Let $j : U \hookrightarrow Y$ be
an open subscheme of $Y$.
Then define a functor
$$
C_U \Phi  \quad : \; CS(X) \lto Sh(X) \qquad \qquad \mbox{as following:}
$$
for any sheaf $\f \in CS(X)$:
$$
C_U \Phi (\f)  \eqdef
 \mathop{\pl}\limits_{k \in {\bf N}}
\Phi (i_* j_* j^* (\f / J^k \f))   \quad \mbox{.}
$$

\begin{nt}
{ \em
It is not difficult to understand that
if the sheaf
$\f$
is a coherent sheaf on
$X$,
then for any
 $k \in {\bf N}$
the sheaf
$i_* j_* j^* (\f / J^k \f)$
is a quasicoherent sheaf on $X$.
In fact,
the sheaf $\f / J^k \f$
is a coherent sheaf on the scheme $X$.
Moreover, the sheaf
$\f / J^k \f$
is a coherent sheaf on the closed subscheme
$(Y, \oo_X/ J^k)$
of the scheme $X$.
Then $j^* (\f/J^k\f)$ is
a coherent sheaf on the scheme
 $(U, \oo_X / J^k \mid_U)$,
$j_* j^* (\f / J^k \f)$ is
a quasicoherent sheaf on the  scheme
$(Y, \oo_X/ J^k)$.
And
since $i_*$
coincides with the direct image functor
from the subscheme $(Y, \oo_X / J^k)$, we see that
$i_*j_* j^* (\f / J^k \f)$ is a quasicoherent
sheaf on $X$.
}
\end{nt}

\smallskip

\begin{lemma}  \label{lem6}
Let $X$ be a noetherian scheme,
$\Phi \: : \: QS(X) \lto Sh(X)$ be
an exact additive functor,
$i: Y \hookrightarrow X$ be a closed
subscheme of the scheme $X$, which is defined by
 the ideal sheaf  $J$ on $X$.
Let $j : U \hookrightarrow Y$ be
an open subscheme of $Y$.
In addition,
suppose the following:
for any point $x \in X$,
for any open $W \subset X$,
$x \in W$,
there exists an affine open subscheme
$V \subset W$, $x \in V$
such that:
\begin{enumerate}
\item           \label{usl1}
 $V \cap U$ is an affine subscheme;
\item            \label{usl2}
for any quasicoherent sheaf
$\f$
on $X$
\begin{equation}  \label{eqn6}
{H}^1 (V, \Phi (\f) ) = 0
\end{equation}
\end{enumerate}
Then
$ C_U \Phi \: : \: CS(X) \lto Sh(X)$
is an exact and additive functor.
\end{lemma}

\begin{nt}
{ \em
For example,
condition~\ref{usl1}
of lemma~\ref{lem6}
is satisfied in the following cases:
\begin{itemize}
\item
$X$ is a  separated scheme,
$U$ is an affine subscheme
(as on a separated scheme the intersection
of two affine open subschemes is an affine subscheme);
\item
$X$ is
a separated scheme,
and $U$ is a complement
to some Cartier divisor in $Y$.
\end{itemize}
}
\end{nt}
{\bf Proof} (of lemma~\ref{lem6}).

Additivity of the functor
$C_U \Phi$
is obvious from the construction.
Let us show exactness.
Let
$$
0 \lto \f_1 \lto \f_2 \lto \f_3    \lto 0
$$
be an exact sequence of coherent sheaves on $X$.

For any $x \in X$,
for any open $W \subset X$
consider an open affine
$V \subset W$,
$V \ni x$
satisfying conditions~\ref{usl1} ¨~\ref{usl2}
of lemma~\ref{lem6}.
Then for the proof of exactness of the functor
$C_U \Phi$
it suffices to show exactness of the following sequence:
\begin{equation}  \label{eqn7}
0
\lto H^0 (V, C_U \Phi (\f_1) )
\lto H^0 (V, C_U \Phi (\f_2) )
\lto H^0 (V, C_U \Phi (\f_3) )
\lto 0    \mbox{.}
\end{equation}
On the other hand,
the following sequence is exact:
$$
0 \lto
\f_1/ J^k \f_2 \cap \f_1
\lto \f_2/ J^k \f_2
\lto \f_3/ J^k \f_3
\lto 0    \mbox{.}
$$
Since the sheaves in the last sequence
are coherent sheaves on the scheme
$(Y, \oo_X / J^k)$,
by lemma~\ref{lem5}
the following sequence is exact:
$$
0 \lto
i_* j_* j^* (\f_1 / J^k \f_2 \cap \f_1)
\lto
i_* j_* j^* (\f_2/J^k \f_2)
\lto
i_* j_* j^* (\f_3/J^k \f_3)
\lto 0                      \mbox{.}
$$
Since we apply the exact functor $\Phi$
to the last sequence,
we obtain exactness of the following sequence from
$Sh(X)$:
$$
0 \lto
\Phi(i_* j_* j^* (\f_1 / J^k \f_2 \cap \f_1))
\lto
\Phi(i_* j_* j^* (\f_2/J^k \f_2))
\lto
\Phi(i_* j_* j^* (\f_3/J^k \f_3))
\lto 0                         \mbox{.}
$$
Now from~(\ref{eqn6})
we obtain exactness of the following sequence:
\begin{flushleft}
$
%\begin{array}{c}
0 \to
H^0(V,\Phi(i_* j_* j^* (\f_1 / J^k \f_2 \cap \f_1)) )
\to
H^0(V,\Phi(i_* j_* j^* (\f_2/J^k \f_2)) )
\to
%\end{array}
$
\end{flushleft}
\begin{equation} \label{eqn8}
\qquad  \qquad \qquad \qquad \qquad \qquad \qquad \qquad \qquad \qquad
\to
H^0(V,\Phi(i_* j_* j^* (\f_3/J^k \f_3)) )
\to 0                         \mbox{.}
\end{equation}
Also for any natural numbers
$k_1 \le k_2$
we have the following exact sequence of sheaves:
$$
0 \lto
J^{k_1} \f_2 \cap \f_1 /
J^{k_2} \f_2  \cap \f_1
\lto
\f_1/
J^{k_2} \f_2  \cap \f_1
\lto
\f_1 /
J^{k_1} \f_2 \cap \f_1
\lto
0  \mbox{.}
$$
Hence, as above,
the sheaf
$
J^{k_1} \f_2 \cap \f_1 /
J^{k_2} \f_2  \cap \f_1
$
is a coherent sheaf on the scheme
$(Y, \oo_X/ J^{k_2})$.
Therefore by lemma~\ref{lem5}
the following sequence is exact:
\begin{flushleft}
$
0 \lto
i_* j_* j^* (J^{k_1} \f_2 \cap \f_1 /
J^{k_2} \f_2  \cap \f_1 )
\lto
i_* j_* j^* (\f_1/
J^{k_2} \f_2  \cap \f_1)
\lto
$
\end{flushleft}
\begin{flushright}
$
\lto
i_* j_* j^* (\f_1 /
J^{k_1} \f_2 \cap \f_1)
\lto
0  \mbox{.}
$
\end{flushright}
Since the functor $\Phi$
is exact, we have  exactness of the sequence:
\begin{flushleft}
$
0 \lto
\Phi(i_* j_* j^* (J^{k_1} \f_2 \cap \f_1 /
J^{k_2} \f_2  \cap \f_1 ) )
\lto
\Phi(i_* j_* j^* (\f_1/
J^{k_2} \f_2  \cap \f_1) )
\lto
$
\end{flushleft}
\begin{flushright}
$
\lto
\Phi(i_* j_* j^* (\f_1 /
J^{k_1} \f_2 \cap \f_1) )
\lto
0  \mbox{.}
$
\end{flushright}
And from~(\ref{eqn6})
we obtain that the following map is a surjective map:
\begin{equation}  \label{eqn9}
\Phi(i_* j_* j^* (\f_1/
J^{k_2} \f_2  \cap \f_1) )
\lto
\Phi(i_* j_* j^* (\f_1 /
J^{k_1} \f_2 \cap \f_1) )
\end{equation}
Now taking the projective limit with respect to all
 $k \in {\bf N}$
and using~(\ref{eqn9}), from which it follows the Mittag-Leffler
condition (see~\cite[ ch.II, \S9]{Ha}),
we obtain  exactness of the following sequence:
\begin{flushleft}
$
0 \to
\mathop{\pl}\limits_{k \in {\bf N}}
H^0(V,\Phi(i_* j_* j^* (\f_1 / J^k \f_2 \cap \f_1)) )
 \to
\mathop{\pl}\limits_{k \in {\bf N}}
H^0(V,\Phi(i_* j_* j^* (\f_2/J^k \f_2)) )
\to
$
\end{flushleft}
\begin{flushright}
$
\to
\mathop{\pl}\limits_{k \in {\bf N}}
H^0(V,\Phi(i_* j_* j^* (\f_3/J^k \f_3)) )
\to 0
$
\end{flushright}
\begin{center}
or
\end{center}
\begin{flushleft}
$
0 \lto
\mathop{\pl}\limits_{k \in {\bf N}}
H^0(V,\Phi(i_* j_* j^* (\f_1 / J^k \f_2 \cap \f_1)) )
\lto
$
\end{flushleft}
\begin{flushright}
$
\lto
H^0(V, C_U \Phi (\f_2) )
\lto
H^0(V, C_U \Phi (\f_3) )
\lto 0  \mbox{.}
$
\end{flushright}
Thus for the proof of  exactness of sequence~(\ref{eqn7})
we have to show that
\begin{equation}  \label{eqn10}
H^0(V, C_U \Phi (\f_1) )
=
\mathop{\pl}\limits_{k \in {\bf N}}
H^0(V,\Phi(i_* j_* j^* (\f_1 / J^k \f_2 \cap \f_1)) ) \mbox{.}
\end{equation}
Since the scheme $X$
is a noetherian scheme,
by Artin-Rees lemma (see~\cite[ cor. 10.10]{AM})
there exists
$l \in {\bf N}$
such that for all
$k \le l$:
\begin{equation}  \label{eqn11}
J^k \f_2 \cap \f_1   =   J^{k-l} (J^l \f_2 \cap \f_1)
\mbox{.}
\end{equation}

From~(\ref{eqn11}) we obtain that
the maps
$$
\f_1 / J^k \f_2 \cap \f_1   \lto
\f_1 / J^{k-l} \f_1
$$
are well defined and surjective.
The same is for
$$
\f_1 / J^k \f_1  \lto
\f_1/ J^k \f_2 \cap \f_1     \quad \mbox{.}
$$
Further, appliing successively  the functors
$j^*$,
$j_*$, $i_*$,
$\Phi$ and $H^0 (V, \cdot)$
and using
lemma~\ref{lem5} and condition~(\ref{eqn6}),
we obtain cofinality of the projective systems:
$$
H^0(V, \Phi(i_* j_* j^* (\f_1 / J^k \f_2 \cap \f_1)) )
\qquad
\mbox{and}
\qquad
H^0(V, \Phi(i_* j_* j^* (\f_1 / J^k \f_1)) )  \mbox{.}
$$
Therefore~(\ref{eqn10})
is satisfied, and lemma~\ref{lem6}
is proved.
%\vspace{0.5cm}

\bigskip

%\vspace{0.3cm}
\begin{lemma} \label{lem7}
Let $X$
be a noetherian scheme,
$\Phi \: : \: QS(X) \lto Sh(X)$ be an
exact additive functor,
$i: Y \hookrightarrow X$ be a closed
subscheme of the scheme $X$, which is defined by
 the ideal sheaf $J$ on $X$.
Let $j : U \hookrightarrow Y$ be
an open subscheme of $Y$ such that
for any point $x \in X$
there exists an affine neighbourhood
$V \ni x$ such that
 $V \cap U$ is an affine subscheme.
Let $\U = \{ U_i \}_{i \in I}  $ be
an affine open covering of the scheme $X$.
Suppose that for any
 $k \ge 1$,
for any coherent sheaf
$\g$ on the scheme $X$  we have
\begin{eqnarray}
\check{H}^m (\U, \, \Phi (i_* j_* j^* (\g/ J^k \g))) & = & 0
\quad  \mbox{for any} \quad m \ge 1 \label{eqn12}        \\
%\end{equation}
%\begin{equation}
H^1 (\bigcap\limits_{i \in I_0} U_i, \,
\Phi (i_* j_* j^* (\g/ J^k \g))) &  = & 0    \quad
\mbox{for any subset}  \quad I_0 \subset I
\label{eqn12a}
      \\
%\end{equation}
%\begin{equation}
H^1(X, \, \Phi (i_* j_* j^* (\g/ J^k \g))) & = &  0   \quad \mbox{.}
\label{eqn12b}
\end{eqnarray}
Then $\check{H}^m (\U, \, C_U \Phi (\f))  = 0 $
for any $m \ge 1$
and for any coherent sheaf
$\f$ on $X$.
\end{lemma}
\proof
Let $\f$
be any coherent sheaf on $X$.
Let
$$
p_n \; : \;
\mathop{\prod\limits_{I_0 \subset I}}\limits_{\sharp I_0 = n+1}
H^0 (\bigcap\limits_{i \in I_0} U_i,  \,
\Phi (i_* j_* j^* (\f/ J^k \f)))
\lto
\mathop{\prod\limits_{I_0 \subset I}}\limits_{\sharp I_0 = n+2}
H^0 (\bigcap\limits_{i \in I_0} U_i,\,
\Phi (i_* j_* j^* (\f/ J^k \f)))
$$
be the map which is arised from the \v{C}ech
complex with respect to the covering
$\U$.
Define $H^n_k \eqdef \Ker p_n$.
In addition,
$$H^0_k = H^0 (X,    \,
\Phi (i_* j_* j^* (\f/ J^k \f)))  \quad \mbox{.} $$
From~(\ref{eqn12})
we obtain at once that for any $n \ge 1$
$$H^n_k = \Image  p_{n-1} \quad \mbox{.}  $$
Therefore for any $m \ge 0$
the following   sequence   is exact:
\begin{equation}  \label{ash}
0  \lto
H^n_k   \lto
\mathop{\prod\limits_{I_0 \subset I}}\limits_{\sharp I_0 = n+1}
H^0 (\bigcap\limits_{i \in I_0} U_i, \,
\Phi (i_* j_* j^* (\f/ J^k \f)))
\stackrel{p_n}{\lto}
H^{n+1}_k
\lto
0
\end{equation}

For any natural numbers $k_1 \le k_2$
we have the exact sequence
$$
0
        \lto
        J^{k_1} \f / J^{k_2} \f
        \lto
        \f/ J^{k_2} \f
        \lto
        \f / J^{k_1} \f
        \lto
        0   \mbox{.}
$$
Since the sheaves of this sequence
are coherent sheaves on the scheme
$(Y, \oo_X / J^{k_2})$,
by lemma~\ref{lem5}
the following sequence is exact:
$$
0 \lto i_* j_* j^*
(       J^{k_1} \f / J^{k_2} \f)
        \lto
        i_* j_* j^*
 (\f/ J^{k_2} \f)
        \lto
        i_* j_* j^*
        (\f / J^{k_1} \f )
        \lto
        0   \mbox{.}
$$
Further, from exactness
of the functor
 $\Phi$ and
condition~(\ref{eqn12b})
we obtain
surjectivity of the following maps for any natural numbers
 $k_1 \le k_2$:
        \begin{equation}  \label{kva}
 H^0(X, \, \Phi( i_* j_* j^*
 (\f/ J^{k_2} \f)))
        \lto
 H^0(X, \,  \Phi( i_* j_* j^*
        (\f / J^{k_1} \f )))  \mbox{.}
         \end{equation}
By condition~(\ref{eqn12a}),
we obtain as well
that for any
$k_1 \le k_2$ and any $I_0 \subset I$
the maps
        \begin{equation}  \label{kva1}
 H^0(\bigcap_{i \in I_0} U_i,   \, \Phi( i_* j_* j^*
 (\f/ J^{k_2} \f)))
        \lto
 H^0(\bigcap_{i \in I_0} U_i,   \, \Phi( i_* j_* j^*
        (\f / J^{k_1} \f )))  \mbox{.}
         \end{equation}
are surjective maps.

Let us prove
that for any
$n \ge 0$   the map
\begin{equation}        \label{sur}
H^n_{k_2}  \lto H^n_{k_1}
\end{equation}
is a surjective map for any
$k_1 \le k_2$.
In the  case $ n = 0$
it is statement~(\ref{kva}).
For arbitrary $n$  it follows
from surjectivity of
$p_{n-1}$
in exact sequence~(\ref{ash})
and, also, surjectivity of map~(\ref{kva1}).

Now taking the projective limit in~(\ref{ash})
and using surjectivity
of~(\ref{sur}),
from which it follows the Mittag-Leffler condition
for the projective systems    $^n_k$,
we obtain
exactness of the following sequence
for any
$n \ge 0$:
$$
        0
        \lto
        \mathop{\pl}\limits_k H^n_k
        \lto
        \mathop{\prod_{I_0 \in I}}\limits_{\sharp I_0 = n+1}
        H^0( \bigcap_{i \in I_0} U_i,
        C_U \Phi(\f))
        \lto
        \mathop{\pl}\limits_k H^{n+1}_k
         \lto 0
          \mbox{.}
        $$
Hence it follows at once that
for any   $m \ge 1$ $\check{H}^m (\U, \, C_U \Phi (\f)) = 0$.
Lemma~\ref{lem7} is proved.

\medskip
\medskip
\medskip
In the sequel we'll use the following variant
of A.~Kartan lemma,
which connects the \v{C}ech cohomologies groups
with the usual  cohomologies groups of sheaves.
For a sheaf $\A$
on a topological space $V$
by $\check{H}^q (V, \A)$
denote direct
 limit of \v{C}ech cohomologies
with respect to all coverings  of the space $V$.

\begin{lemma}  \label{Kartan}
Let $X$
be a topological space
and $\A$
be a sheaf on $•$.
Suppoce that it is possible to cover $X$
by family {$\bf U$} of open sets
such that this family has the following properties:
\begin{enumerate}
\item If $\bf U$
contains $U'$ and $U''$,
then it contains $U' \cap U''$;
\item  $\bf U$
contains
arbitrarily
 small open sets;
\item  $\check{H}^q (U, \A) = 0$
for any $q \ge 1$  and $U \in {\bf U}$.
\end{enumerate}
Under these conditions we have isomorphism:
$$
\check{H}^q (X, \A) \lto H^q(X, \A)  \mbox{.}
$$
\end{lemma}
\proof See theorem~5.9.2 of~\cite[ ch.2]{G}.

\bigskip
\begin{lemma}  \label{lem8}
Let $X$
be a noetherian scheme,
$j : U \hookrightarrow X$
be an open affine subscheme.
Then for any quasicoherent sheaf
 $\f$  on $U$:
\begin{eqnarray}
 H^0(X, j_* \f) & =  & H^0 (U, \f)  \label{eqn17}  \\
 H^i (X, j_* \f) & = &  0  \qquad \qquad
\mbox{if} \quad i > 0   \mbox{.}
 \label{eqn18}
\end{eqnarray}
 \end{lemma}
\proof

Equality~(\ref{eqn17})
follows from the construction of the functor
$j_*$.
Let us prove~(\ref{eqn18}).
Embed the quasicoherent sheaf
$\f$
in a flasque quasicoherent sheaf
 $\g$ on $U$.
(It always can do, see~\cite[ ch.~III, \S3]{Ha}.)
$$
0
\lto
\f
\lto
\g
\lto
\g  /  \f
\lto 0
$$
Now the following sequence is exact:
\begin{equation}  \label{eqn19}
0 \lto j_* \f
\lto
j_* \g
\lto
j_* (\g  /\f)
\lto
0           \mbox{.}
\end{equation}
(Indeed,
$R^1 j_* \f =0$.
The last follows from quasicoherentness of the sheaf $R^1 j_* \f$
and for any affine open
$V \subset X$:
$H^0 (V, R^1 j_* \f) =
H^1 (V \cap U, \f) = 0
$, as  from
separateness of
$X$ it follows that
$V \cap U$ is an affine scheme.)
Besides,
it is not difficult to see that the sheaf
$j_* \g$ is an flasque sheaf.
Therefore
\begin{equation}   \label{eqn20}
H^i (X, j_* \g)= 0  \qquad \qquad
\mbox{for any} \quad i > 0 \mbox{.}
\end{equation}
Besides,
the map
$$
H^0(X, j_* \g)  \lto H^0 (X, j_* (\g / \f))
$$
is surjective,
as $H^0 (X, j_* \g)  =
H^0 (U, \g)$,
$H^0 (X, j_* (\g / \f))
= H^0 (U, \g / \f)  \mbox{.}
$
And from afinneness of $U$
it follows that
$H^1 (U, \f) = 0$,
therefore the following map is surjective:
$$
H^0 (U, \g) \lto H^0 (U, \g / \f) \mbox{.}
$$
Hence and from~(\ref{eqn20})  we obtain that
$$
H^1 (X, j_* \f) = 0   \mbox{.}
$$
Further, if $i > 1$,
then from~(\ref{eqn20})
and from the long cohomological sequence associated with~(\ref{eqn19})
it follows that
$$
H^i(X, j_* \f) = H^{i-1}(X, j_* (\g / \f))  \mbox{.}
$$
But the sheaf $\g / \f$  is a quasicoherent sheaf
on  $U$.
Hence, by induction, it is possible to assume that
$$
H^{i-1} (X, j_* (\g / \f)) = 0  \mbox{.}
$$
Therefore $H^i(X, j_* \f) = 0$.
Lemma~\ref{lem8} is proved.

\bigskip

\begin{lemma}       \label{lem9}
Let $X$ be a noetherian scheme.
Let
$Y \hookrightarrow X$ be
a closed subscheme,
which is defined by the ideal sheaf $J$,
and $j : U \hookrightarrow X$ be the
open subscheme
which is complement to the subscheme $Y$.
Let a sheaf $\f$ be a quasicoherent sheaf on $X$, and
consider the following exact sequence of
quasicoherent sheaves on $X$:
$$
0 \lto \h \lto \f \lto j_* j^* \f
\lto \g \lto 0  \mbox{}
$$
which is induced by the natural map
$
\f \lto j_* j^* \f
$.

Let  $\h = \mathop{\il}\limits_i  \h_i
  $     and  $
\g = \mathop{\il}\limits_i  \g_i
  $,
$i \in I
$,
where the sheaves
 $\h_i$ ¨ $\g_i$
are coherent sheaves on  the scheme
 $X$
for any $i \in I$.

Then for any $i \in I$
there exists $l(i) \in {\bf N}$ such that
\begin{equation}      \label{hhhh}
J^{l(i)} \cdot \h_i = 0
\qquad   \mbox{and}  \qquad
J^{l(i)} \cdot \g_i = 0  \mbox{.}
\end{equation}
\end{lemma}
\proof
From the exact sequence
$$
0 \lto
J^m
\lto
\oo_X
\lto   \oo_X/J^m   \lto 0
$$
it follows the following sequence of quasicoherent sheaves on $X$
$$
0 \lto
{\cal H}om_{X} (\oo_X/J^m, \f)
\lto
{\cal H}om_{X} (\oo_X, \f)
\lto
{\cal H}om_{X}  (J^m, \f)
\lto
{\cal E}xt^1_{X}
(\oo_X/J^m, \f)  \mbox{.}
$$
Taking direct limit with respect to $m$, we obtain
$$
0 \lto
\mathop{\il}_m \,
{\cal H}om_X (\oo_X/J^m, \f)
\lto
\f
\lto
\mathop{\il}_m   \,
{\cal H}om_X (J^m, \f)
\lto
\mathop{\il}_m  \,
{\cal E}xt^1_X
(\oo_X/J^m, \f)
\mbox{.}
$$
By
\cite[ ch.~III, ex.~3.7(a)]{Ha} we have
$$
\mathop{\il}_m   \,
{\cal H}om_X (J^m, \f)
=
j_* j^* \f  \mbox{.}
$$
Now (\ref{hhhh})
follows from
$$
{\cal H} =
\mathop{\il}_m  \,
{\cal H}om_X (\oo_X/J^m, \f)
\qquad
\mbox{and}
\qquad
\g \hookrightarrow
\mathop{\il}_m   \,
{\cal E}xt^1_X
(\oo_X/J^m, \f)     \mbox{.}
$$
Lemma~\ref{lem9} is proved.

\bigskip

%***************************************************************
\section{Construction and its original properties}
Let $X$ be
a noetherian separated scheme.
Consider a flag of closed subschemes
$$
X \supset Y_0 \supset Y_1 \supset \ldots \supset Y_n
$$
in the scheme $X$.
Let $J_j$ be the ideal sheaf of the subscheme $Y_j$ in $X$ ($0 \le j \le n$).
Let $i_j$ be the embedding of the subscheme
 $Y_j  \hookrightarrow X$.
Let $U_i$ be an open subscheme of  $Y_i$
which is complement to the closed subscheme $Y_{i+1}$
 ($0 \le i \le n-1$).
Let $j_i : U_i \hookrightarrow Y_i$ be the open embedding
of the subscheme $U_i$
to the scheme $Y_i$ ($0 \le i \le n-1$).
By definition, let
 $U_n = Y_n$
and $j_n$ be
the identity morphism from
 $U_n$  to $Y_n$.

Assume that for any point
 $x \in X$
there exists an open affine neithbourhood
 $U \ni x$
such that
$U \cap U_i$ is an affine scheme for any $0 \le i \le n$.
In the sequel we'll say that a flag of subschemes
 $\{Y_i, \; 0 \le i \le n \}$ with such condition
is the flag with { \em locally affine complements.}

\begin{nt}
{ \em
For example, the last condition of locally affineness of complements
is satisfied in the following cases
\begin{itemize}
\item $Y_{i+1}$
is the Cartier divisor on the scheme
$Y_i$ ($0 \le i \le n-1$), or
\item
$U_i$ is an affine scheme for any  $0 \le i \le n-1$.
(On a separated scheme the intersection of two open affine subschemes
is an  affine subscheme.)
\end{itemize}
}
\end{nt}

\bigskip
Consider the $n$-dimensional
simplex
and its standard simplicial set (without degenerations).
To be precise, consider the set:
$$
(\{ 0\}, \{ 1\}, \ldots, \{ n\})    \mbox{.}
$$
(Here are all the integers between $0$ and $n$.) \\
Then the simplicial set $S = \{ S_k \}$ :
\begin{itemize}
\item $S_0 \eqdef \{\eta \in \{ 0\}, \{ 1\}, \ldots, \{ n\}   \} $.
\item
$ S_k \eqdef \{ (\eta_0, \ldots, \eta_k),
\quad \mbox{where} \quad \eta_l \in S_0  \quad \mbox{and} \quad
\eta_{l-1} < \eta_l  \} $.
\end{itemize}
The boundary map $\partial_i$ ($0 < i < k$)
is given by eliminating the
$i$-th component of the vector  $(\eta_0, \ldots, \eta_k)$.
(It is the  $i$-th face  of  $(\eta_0, \ldots, \eta_k)$.)
\vspace{0.5cm}  \\
{\bf Definition.}
{\em
For any  $(\eta_0, \ldots , \eta_k)  \in S_k$
define the functor
$$
V_{(\eta_0, \ldots, \eta_k)} \; : \; QS(X)  \lto Sh(X)
$$
uniquely determined by the following inductive conditions:
\begin{enumerate}
\item
$V_{(\eta_0, \ldots, \eta_k)}  $
commutes with direct limits.
\item
If $\f$ is a coherent sheaf, and $\eta \in S_0$,
then
$$
V_{\eta}(\f) \eqdef
\mathop{\pl}\limits_{m \in \bf{N}}
(i_{\eta})_* (j_{\eta})_*
(j_{\eta})^* (\f / J^m_{\eta} \f) \mbox{.}
$$
\item
If $\f$ is a coherent sheaf,
and $(\eta_0, \ldots, \eta_k) \in S_k$ ($k \ge 1$),
then
$$
V_{(\eta_0, \eta_1, \ldots, \eta_k)}(\f) \eqdef
\mathop{\pl}\limits_{m \in {\bf N}}
V_{(\eta_1, \ldots, \eta_k)}
\left( (i_{\eta_0})_*  (j_{\eta_0})_*  (j_{\eta_0})^*
(\f / J^m_{\eta_0} \f) \right)  \mbox{.}
$$
\end{enumerate}
}
In the sequel, to avoid the confusion of notations
in the case of
a lot of schemes and flags of closed subschemes
 we'll use sometimes
the equivalent  notation for
$V_{(\eta_0, \ldots, \eta_k)}(\f)$,
in which the closed subschemes is written explicitly:
$$
V_{(\eta_0, \ldots, \eta_k)}(\f) =
V_{(Y_{\eta_0}, \ldots, Y_{\eta_k})}(X, \f) \mbox{.}
$$

\begin{prop} \label{predl1}
Let $ \sigma = (\eta_0, \ldots, \eta_k) \in S_k$.
Then
\begin{enumerate}
\item       \label{pun1}
The functor $V_{\sigma} :
QS(X) \lto Sh(X)$ is well defined.
\item        \label{pun2}
The functor $V_{\sigma}$
is exact and additive.
\item     \label{pun3}
The functor $V_{\sigma}$
is local on $X$, i.~e.,
for any open
$U \subset X$
for any quasicoherent sheaf
 $\f$ on $X$:
$$
V_{(Y_{\eta_0}, \ldots, Y_{\eta_k})}(X, \f)  \mid_U =
V_{(Y_{\eta_0} \cap U, \ldots, Y_{\eta_k} \cap U  )} (U, \f \mid_U)  \mbox{.}
$$
(Here if  $Y_j \cap U = \o$,
then $Y_i \cap U$ is an empty subscheme of $U$ which is
defined by the ideal sheaf $\oo_U$.)
\item For any quasicoherent sheaf $\f $ on the scheme $X$
the sheaf $V_{(\eta_0, \ldots, \eta_k)}(\f)$ is a sheaf of
 $\oo_X$-modules
with the support  on the subscheme $Y_{\eta_k}$.
(Usually, this sheaf  is not quasicoherent.)
\item   \label{pun5}
For any quasicoherent sheaf
 $\f$ on $X$:
$$
V_{\sigma}(\f) =
V_{\sigma}(\oo_X)  \otimes_{\oo_X} \f \mbox{.}
$$
\item     \label{pun6}
If all
 $U_i$ is affine ($0 \le i \le n$),
then for any affine covering
 $\U$  of the scheme $X$,
for any quasicoherent sheaf
 $\f$ on $X$,
for any $m \ge 1$:
$$
\check{H}^m (\U, V_{\sigma}(\f)) ) = 0 \mbox{.}
$$
\item    \label{pun7}
If all $U_i$ is affine ($0 \le i \le n$),
then for any quasicoherent
sheaf $\f$ on $X$, for any $m \ge 1$
$$
H^m (X, V_{\sigma}(\f)) = 0    \mbox{.}
$$
\end{enumerate}
\end{prop}
\proof\\
1.
Well-posedness of the definition of
$V_{\sigma}$  is proved by induction
by means of using of lemma~\ref{lem1},
lemma~\ref{lem2}, lemma~\ref{lem3},
lemma~\ref{lem6}, lemma~\ref{lem7} and lemma~\ref{Kartan}.
Let us check the base of induction for lemmas~\ref{lem6}  and~\ref{lem7}
(when the functor $\Phi = id $). Namely
\begin{itemize}
\item[a)]
for any affine scheme
$\V$
and any quasicoherent sheaf
 $\f$ on $V$
$$
H^1(V, \f) = 0 \mbox{.}
$$
\item[b)]
Let us show
that if
 $i : Y \hookrightarrow X$
is a closed subscheme with
the ideal sheaf
 $J$,
$j : U \hookrightarrow Y$ is an open imbedding
of the affine scheme
  $U$ in $Y$.
Then for any quasicoherent sheaf
 $\f$
on $X$,
for any $k \ge 1$,
for any affine open covering $\U$
of the scheme $X$, for any $m \ge 1$:
$$
\check{H}^m (\U, i_* j_* j^* (\f / J^k \f))= 0  \mbox{.}
$$
From affineness of the covering
$\U$ it follows that it is acyclic for quasicoherent sheaves.
Consequently the \v{C}ech cohomologies groups with respect to this covering
coincide with the usual cohomologies groups of
quasicoherent sheaves. (See~\cite[ ch.3, theorem~4.5]{Ha})
Therefore it suffices to prove that
for any integer
 $k \ge 1$
\begin{equation} \label{dada}
H^m (X, i_* j_* j^* (\f / J^k \f)) = 0  \mbox{.}
\end{equation}
If $k =1$,
then the sheaf $\f / J\f$
is quasicoherent with respect to the subscheme
$Y$, and
$$
H^m (X, i_* j_* j^* (\f / J \f))=
H^m (Y, j_* (j^* (\f / J \f))) = 0 \mbox{.}
$$
Where the last equality follows from lemma~\ref{lem8}.

If $k > 1$,
then by lemma~\ref{lem5}
the following sequence is exact
\begin{equation}       \label{da}
0 \lto  i_* j_* j^* (J^{k-1} \f / J^k \f)
\lto
i_* j_* j^* (\f / J^k \f)
\lto
i_* j_* j^* (\f / J^{k-1} \f)
\lto
0  \mbox{.}
\end{equation}
In addition, the sheaves $J^{k-1} \f / J^k \f$
and $\f / J^{k-1} \f$
are quasicoherent with respect to the subscheme
 $(Y, \oo_X / J^{k-1})$.
Therefore, by induction,  we obtain
$$
H^m(X, i_* j_* j^* (J^{k-1} \f / J^k \f)) = 0  \qquad  \quad
\mbox{and}
\qquad \quad
H^m(X, i_* j_* j^* (\f/ J^{k-1} \f)) =0  \mbox{.}
$$
Hence, from~(\ref{da})  we have~(\ref{dada}).
Item~\ref{pun1} of proposition~\ref{predl1} is proved.
\end{itemize}
2.
The proof of this item is analogous to the proof of item~\ref{pun1}
by means of the same lemmas.\\[0.3cm]
3.
Localness follows by induction from the construction
of the functor
$V_{(\eta_0, \ldots, \eta_k)}$. \\[0.3cm]
4.
This item follows by induction from the construction. \\[0.3cm]
5.
We have the natural map:
$$
\f \lto  V_{\sigma}(\f)  \mbox{,}
$$
which induces the following map:
%\begin{flushleft}
%$
\begin{equation}  \label{dp}
V_{\sigma}(\oo_X)  \otimes_{\oo_X} \f
\lto
V_{\sigma}(\oo_X)  \otimes_{\oo_X}
V_{\sigma}(\f)
\lto
%$
%\end{flushleft}
%\begin{equation}  \label{dp}
%\qquad \qquad \qquad
%\lto
V_{\sigma}(\oo_X)  \otimes_{
V_{\sigma}(\oo_X)
}
V_{\sigma}(\f)
=
V_{\sigma}(\f)  \mbox{.}
\end{equation}
\vspace{0.3pt}
Let us show that~(\ref{dp})
gives us an isomorphism between
$V_{\sigma} (\oo_X) \otimes_{\oo_X} \f$   and
$V_{\sigma} (\f)$.
Since the functor $V_{\sigma} $
and tensor products commute with direct limits,
we can assume that $\f$ is a coherent sheaf.
In view of item~\ref{pun3} of this proposition,
we can restrict ourself to the local situation.
That is,  we suppose
 $X = \Spec A$,
$\f = \tilde{M}$
for some finitely generated
$A$-module $M$.
Then for some $r \in \bf{N}$
there exists an exact sequence of sheaves
as:
$$
0
\lto
\tilde{N}
\lto
\oo_X^{\oplus r}
\lto
\tilde{M}
\lto
0   \mbox{,}
$$
where $N$ is some finitely generated
$A$-module.
Hence we obtain the commutative diagram:
$$
\begin{array}{ccccccc}
&
V_{\sigma} (\oo_X)
\otimes_{\oo_X} \tilde{N}
 & \lto
&
V_{\sigma} (\oo_X)
\otimes_{\oo_X}
\oo_X^{\oplus r} &
\lto &
V_{\sigma} (\oo_X)
\otimes_{\oo_X}
\tilde{M}  &
\lto
0
\\           &
 \begin{picture}(0,28)
 \put(0,26){\vector(0,-1){26}}
 \put(0,11){$\; \gamma $}
 \end{picture}
 &&
 \begin{picture}(0,28)
 \put(0,26){\vector(0,-1){26}}
 \put(0,10){$\; \beta $}
 \end{picture}
 &&
 \begin{picture}(0,28)
 \put(0,26){\vector(0,-1){26}}
 \put(0,11){$\; \alpha $}
 \end{picture}
         \\
 0 \lto
 &  V_{\sigma} (\tilde{N})
     &
\lto &
V_{\sigma} (\oo_X^{\oplus r})
&
\stackrel{\delta}{\lto}
&
V_{\sigma}(\tilde{M})
&
\lto 0 \mbox{,}
\end{array}
$$
where the lower row is exact by virtue of item~\ref{pun2}.
Besides, it is clear that $\beta$
is an isomorphism.
Therefore from surjectivity
of $\delta$
it follows that $\alpha$ is surjective.
Since $\tilde{N}$ is a coherent sheaf,
we have that the map $\gamma$
is surjective as well.
Hence, from exactness of the lower row
and non complicated diagram search
it follows that the map $\alpha$ is injective.  \\[0.3cm]
6.
The proof is similar to
the proof of item~\ref{pun1}
by means of the same lemmas.\\[0.3cm]
7.
This item follows from the previous item of this proposition
and lemma~\ref{Kartan}.
(Since every point has arbitrarily small
affine neithbourhood  with affine intersection to
all $U_i$,
we obtain that this affine neithbourhood satisfies
item~\ref{pun6} of proposition~\ref{predl1}.)
\vspace{0.5cm}
\begin{prop}  \label{predl2}
\begin{enumerate}
\item
Let $X$
be a noetherian separated scheme.
Let
$$Y_0 \supset Y_1 \supset \ldots \supset Y_n
\qquad \qquad \mbox{and} \qquad \qquad
Y'_0 \supset Y'_1 \supset \ldots \supset Y'_n $$
be two flags of closed subschemes in $X$
with the corresponding ideal sheaves
$J_i$  and $J'_i$ ($0 \le i \le n$)
such that for any
 $0 \le i \le n$
there exist integers
$l_i \ge 1$
and $l'_i \ge 1$ with the following properties:
\begin{equation} \label{eqn23}
J_i^{l_i} \subset J'_i   \qquad \qquad \mbox{and} \qquad \qquad
(J'_i)^{l'_i} \subset J_i  \mbox{.}
\end{equation}
Then the functors
$$
V_{(Y_{\eta_0}, \ldots, Y_{\eta_{k}})}(X, \cdot)
\qquad \qquad
\mbox{and}
\qquad \qquad
V_{(Y'_{\eta_0}, \ldots, Y'_{\eta_{k}})}(X, \cdot)
$$
coincides for any
 $(\eta_0, \ldots, \eta_k) \in S_k$.
\item
Consider the flag of closed subschemes:
$$
X \supset  Z \supset Y_0 \supset \ldots \supset Y_n
$$
on a noetherian separated scheme $X$.
Let $i: Z \hookrightarrow X$ be a closed imbedding.
Then for any quasicoherent sheaf $\f$
on the scheme $Z$ we have
$$
i_*  \left( V_{(Y_0, \ldots, Y_n)} (Z, \f) \right)  =
V_{(Y_0, \ldots, Y_n)} (X, i_* \f) \mbox{.}
$$
\end{enumerate}
\end{prop}

\begin{nt}  {\em
Condition~(\ref{eqn23})
is equivalent to the statement that
the topological spaces of subschemes
 $Y_i$ and $Y'_i$
are the same for all $0 \le i \le n$.
} \end{nt}
{\bf Proof} (of proposition~\ref{predl2}).\\
1.
It suffices to prove for the case when
$J_i = J'_i$ for all $i \ne j$, $0 \le i \le n$,
where some fixing $0 \le j \le n$.
From inductance of the definition of the functor
$V_{(Y_{\eta_0}, \ldots, Y_{\eta_k})} (X, \cdot)$
 (and $V_{(Y'_{\eta_0}, \ldots, Y'_{\eta_k})} (X, \cdot)$)
we can restrict ourself to the case  $j = 0$.
This case follows from cofinality of the  projective systems
$\f/ J_j^k$  and $\f / (J'_j)^k$ (from condition~\ref{eqn23}).\\[0.3cm]
2.
This item follows at once from the construction and the fact that
supports of all appearing from induction sheaves
are on the subscheme $Z$.

%*************************************************************
\section{Complexes and their exactness}
Consider again the usual $n$-simplex without degenerations
$S= \{ S_k, \; 0 \le k \le n \}$.
If $\sigma= (\eta_0, \ldots, \eta_k) \in S_k$,
then $\partial_i (\sigma)$ is the $i$-th face of $\sigma$ ($0 \le i \le k$).
Then  {\em define} the morphism of functors
$$
d_i(\sigma) \quad : \quad  V_{\partial_i(\sigma)} \lto V_{\sigma} \mbox{,}
\qquad \quad \mbox{as}
$$
commuting with direct limits, and
on coherent sheaves it is a map
\begin{equation}   \label{tank}
V_{\partial_i(\sigma)}(\f) \lto V_{\sigma}(\f)  \mbox{}
\end{equation}
which is defined by the following rules:
\begin{itemize}
\item[a)]
if $i =0$,
then (\ref{tank})
is obtained from application of the functor
$V_{\partial_0 (\sigma)}$
to the map
$$
\f \lto
(i_{\eta_0})_* (j_{\eta_0})_*  (j_{\eta_0})^*
(\f / J_{\eta_0}^m \f)
$$
and passage to the projective limit on $m$;
\item[b)]
if $i=1$, $k=1$,
then we have the natural map
$$
(i_{\eta_0})_* (j_{\eta_0})_*  (j_{\eta_0})^*
(\f / J_{\eta_0}^m \f)
\lto
V_{(\eta_1)} ((i_{\eta_0})_* (j_{\eta_0})_*  (j_{\eta_0})^*
(\f / J_{\eta_0}^m \f)) \mbox{.}
$$
Now after passage to the projective limit on
 $m$
we obtain the map~(\ref{tank}) in this case.
\item[c)]
$i \ne 0$, $k >1$,
then  from induction on
$k$ we can suppose that we have the map
$$
V_{\partial_{i-1} \cdot
(\partial_0 (\sigma))
} ((i_{\eta_0})_* (j_{\eta_0})_*  (j_{\eta_0})^*
(\f / J_{\eta_0}^m \f))
\lto
V_{
\partial_0 (\sigma)
} ((i_{\eta_0})_* (j_{\eta_0})_*  (j_{\eta_0})^*
(\f / J_{\eta_0}^m \f))  \mbox{.}
$$
And passage to the projective limit on $m$
gives us the map~(\ref{tank})
in this case.
\end{itemize}

\begin{prop}     \label{predl3}
For any $1 \le k \le n$, $0 \le i \le k$ define
$$
d_i^k  \eqdef \sum_{\sigma \in S_k} d_i(\sigma)
\quad : \quad \bigoplus_{\sigma \in S_{k-1}} V_{\sigma}
\lto
\bigoplus_{\sigma \in S_k} V_{\sigma}  \mbox{.}
$$
Also define
$$d_0^0 \quad : \quad id \lto \bigoplus_{\sigma \in S_0}
V_{\sigma} \mbox{}$$
as the direct sum of the natural maps
$\f \lto
V_{\sigma} (\f)$. (Here $id$ is the functor
of the natural imbedding of
 $QS(X)$ into $Sh(X)$,
$\f$ is a quasicoherent sheaf on $X$,
$\sigma \in S_0$.)

Then for all $0 \le i < j \le k \le n-1$ we have
\begin{equation}  \label{npp}
 d_j^{k+1} d_i^k  = d_i^{k+1} d_{j-1}^k   \mbox{.}
\end{equation}
\end{prop}
{\bf Proof}.
Using the inductance of the definition, the proof is done by induction from
non complicated consideration
of some cases.
It suffices to consider the small $i$ and $k$ only.
(For example, see
similar cases in~\cite[ \S 2.4]{H2}
or~\cite{H1}.) \\[0.5cm]

{\em Define}
$$
d^m \eqdef \sum\limits_{0 \le i \le m} (-1)^i d^m_i
$$
Then proposition~\ref{predl3}  makes possible
to construct
the complex of sheaves  $V(\f)$  from any quasicoherent sheaf
 $\f$ on $X$ in the  standard way:
$$
\ldots
\lto
\bigoplus_{\sigma \in S_{m-1}} V_{\sigma}(\f)
\stackrel{d^m}{\lto}
\bigoplus_{\sigma \in S_m} V_{\sigma}(\f)
\lto
\ldots
$$
Where  $d^{m+1} d^m = 0$ follows
from~(\ref{npp})
by means of non complicated  direct calculations.

\medskip
\begin{th}   \label{teorem1}
Let $X$
be a noetherian separated scheme. Let
$Y_0 \supset Y_1 \supset \ldots  \supset  Y_n$ be a flag
of closed subschemes with locally affine complements.
Assumee that $Y_0 = X$.
Then the following complex is exact:
\begin{equation}  \label{kff}
0 \lto \f \stackrel{d^0}{\lto} V(\f) \lto 0     \mbox{.}
\end{equation}
\end{th}
\proof
It suffices to consider only the case when the sheaf $\f$ is coherent.
Consider the exact sequence of sheaves
\begin{equation}  \label{krya}
0 \lto \h \lto \f \lto (j_0)_* (j_0)^* \f  \lto \g \lto 0
\end{equation}
Here  $\h$  and $\g$ is
the kernel and the cokernel of
the natural map of sheaves $\f \lto (j_0)_* (j^0)^* \f $.
From exactness of functors $V_{\sigma}$
(for any $\sigma$)
we obtain the following exact sequence of complexes of sheaves:
\begin{equation} \label{kp}
0 \lto V(\h) \lto V(\f) \lto V((j_0)_* (j_0)^* \f)  \lto V(\g) \lto 0
\end{equation}
By lemma~\ref{lem9}
the supports of sheaves $\h$ and $\g$
are on $Y_1$, therefore
in the case $\eta_0 = 0$ we have
$V_{(Y_{\eta_0}, \ldots, Y_{\eta_k})} (X, \h) = 0$ and
$V_{(Y_{\eta_0}, \ldots, Y_{\eta_k})} (X, \g) = 0$.
 Therefore, using it, lemma~\ref{lem9} (which decompose
the sheaves $\h$ and $\g$
in direct limits of sheaves which is coherent on
subschemes with topological space $Y_1$),
proposition~\ref{predl2}, permutability
of the functors $V_{\sigma}$
with direct limits,
we can apply induction on the length of flag  and suppose  that
the complexes
\begin{equation} \label{kh}
0 \lto \h \stackrel{d^0}{\lto} V(\h) \lto 0     \quad  \mbox{and}
\end{equation}
\begin{equation}      \label{kg}
0 \lto \g \stackrel{d^0}{\lto} V(\g) \lto 0
\end{equation}
are already exact.
It is not difficult to understand
that for any
$\sigma = (\eta_0, \ldots, \eta_k) \in S_k$
if $\eta_0 = 0$, then $V_{\sigma}((j_0)_* (j_0)^* \f) =
V_{\sigma}(\f)$;
if $\eta_0 \ne 0$, then
 $V_{\sigma}((j_0)_* (j_0)^* \f) =
V_{\sigma'}(\f)$, where
$\sigma' = (0, \eta_0, \ldots, \eta_k) \in S_{k+1}$.
Hence,
the complex
$
V((j_0)_* (j_0)^* \f)
$
has the same components   $V_{\sigma'}(\f)$
in the degree
$k$ and $k+1$.
Therefore, successively from the highest degrees spliting off
the trivial complexes
$$
0 \lto  V_{\sigma}((j_0)_* (j_0)^* \f)
\lto  V_{\sigma'}((j_0)_* (j_0)^* \f)    \lto 0    \mbox{,}
$$
we obtain exactness of the complex
\begin{equation}   \label{kf}
0 \lto
(j_0)_* (j_0)^* \f
\stackrel{d^0}{\lto} V((j_0)_* (j_0)^* \f) \lto 0   \mbox{.}
\end{equation}
Now, since complexes~(\ref{kh}), (\ref{kg}),
(\ref{kf}) are exact,
we obtain exactness of  complex~(\ref{kff})
from exactness of~(\ref{kp}) and (\ref{krya}).
Theorem~\ref{teorem1} is proved.

\medskip
\bigskip
\bigskip

For any $\sigma \in S_k$ define
$$
A_{\sigma}(\f) \eqdef  H^0(X, V_{\sigma} (\f))  \mbox{.}
$$
\begin{prop} \label{predl4}
Let $X$ be a noetherian separated scheme.
$Y_0 \supset Y_1 \supset \ldots  \supset  Y_n$ be a flag
of closed subschemes
such that all
$U_i$ are affine ($0 \le i \le n$).
Let $\sigma \in S_k$ be arbitrary.
Then
\begin{enumerate}
\item        \label{punkt1}
$A_{\sigma}$ is an exact and additive
functor: $QS(X) \lto Ab$.
\item  If $X = \Spec A$,
$M$ is some $A$-module,
then
$$
A_{\sigma}(\tilde{M}) = A_{\sigma} (\oo_X)  \otimes_A M  \mbox{.}
$$
\end{enumerate}
\end{prop}
\proof \\
1.
This item follows at once from items~\ref{pun2} and~\ref{pun7}
of proposition~\ref{predl1}.\\[0.3cm]
2.
Similarly to the proof of item~\ref{pun5}  of proposition~\ref{predl1}
we can suppose that the module $M$
is finitely generated over $A$. Now consider
the exact sequence of  $A$-modules:
$$
0 \lto N \lto A^{\oplus^r} \lto M \lto 0  \mbox{.}
$$

Hence we obtain the commutative diagramm:
$$
\begin{array}{ccccccc}
&
A_{\sigma}(\oo_X)
\otimes_{A} \tilde{N}
 & \lto
&
A_{\sigma}(\oo_X)^{\oplus^r}
 &
\lto &
A_{\sigma} (\oo_X) \otimes_A M
  &
\lto
0
\\           &
 \begin{picture}(0,28)
 \put(0,26){\vector(0,-1){26}}
 \put(0,11){$\; \gamma $}
 \end{picture}
 &&
 \begin{picture}(0,28)
 \put(0,26){\vector(0,-1){26}}
 \put(0,10){$\; \beta $}
 \end{picture}
 &&
 \begin{picture}(0,28)
 \put(0,26){\vector(0,-1){26}}
 \put(0,11){$\; \alpha $}
 \end{picture}
         \\
 0 \lto
 & A_{\sigma} (\tilde{N})
     &
\lto &
A_{\sigma} (\oo_X^{\oplus r})
&
\stackrel{\delta}{\lto}
&
A_{\sigma}(\tilde{M})
&
\lto 0 \mbox{,}
\end{array}
$$
where the lower row is exact by virtue of item~\ref{punkt1}
of this proposition.
It is clear that $\beta$ is an isomorphism.
Therefore, arguing as in item~\ref{pun5},
we obtain at first surjectivity of the map~$\alpha$,
and afterwards we obtain injectivity of~$\alpha$.
 Proposition~\ref{predl4} is proved.
%\vspace{0.5cm}

\medskip
\bigskip
\bigskip
Let $\f$ be any quasicoherent sheaf on  $X$.
Apply the functor
 $H^0(X, \cdot)$ to the complex $V(\f)$.
We obtain the complex of abelian groups $A(\f)$:
$$
\ldots \lto \bigoplus_{\sigma \in S_{m-1}} A_{\sigma}(\f)
\lto
\bigoplus_{\sigma \in S_{m}} A_{\sigma}(\f)
\lto
\ldots   \mbox{.}
$$

\medskip
\begin{th}    \label{teorem2}
Let $X$ be a noetherian separated scheme.  Let
$Y_0 \supset Y_1 \supset \ldots \supset Y_n$ be a flag
of closed subschemes
such that
 $Y_0 = X$ and
all $U_i$ are affine   ($0 \le i \le n$).
Then cohomology of the  complex
 $A(\f)$
coincide with cohomology of the  sheaf
$\f$  on $X$,
i.~e.,  for any $i$
$$
H^i(X, \f) = H^i(A(\f))   \mbox{.}
$$
\end{th}
\proof
From theorem~\ref{teorem1}
and item~\ref{pun7}  of proposition~\ref{predl1}
it follows that $V(\f)$ is
an acyclic resolution for the sheaf
$\f$.
Therefore
it is possible to calculate
cohomology   of the sheaf
 $\f$
by means of global sections of this resolution.
Theorem~\ref{teorem2} is proved. \\
\medskip \\
From the last theorem we obtain at once the following
geometrical corollary.
\begin{th}     \label{teorem3}
Let $X$ be a projective algebraic scheme
of dimension $n$ over a field.
Let
$Y_0 \supset Y_1 \supset \ldots \supset Y_n$ be
a flag of closed subschemes
such that
 $Y_0 = X$ and
$Y_i$
is an ample divisor on the scheme
$Y_{i-1}$ for any
$1 \le i \le n$.
Then for any quasicoherent sheaf
 $\f$
on $X$,
 for any $i$  we have
$$
H^i(X, \f) = H^i(A(\f))   \mbox{.}
$$
\end{th}
\proof In fact, since
$Y_i$ is an ample divisor on $Y_{i-1}$
for all $1 \le i \le n$,we have that
 $U_i$ is an affine scheme for all $0 \le i \le n-1$.
Since $ \dm  Y_n =0$, we have that $U_n= Y_n$ is affine as well.
Now application of theorem~\ref{teorem2}
concludes the proof.
%\bigskip
%\bigskip                                            \\

\smallskip
\medskip
\begin{nt}
{\em
Let us remark that for any quasicoherent sheaf
 $\f$,
for any $\sigma = (\eta_0) \in S_0$     \quad
$A_{\sigma} (\f)$
is the group of section over $U_{\eta_0}$
of the sheaf $\f$
lifted to the formal neighbourhood
of the subscheme $Y_{\eta_0}$
in $X$.
And the complex
 $A(\f)$ can be  interpreted as
the \v{C}ech complex for the  such ''covering ''
of the scheme $X$.
}
\end{nt}

%******************************************************
\section{Combinatorial properties and the Krichever map.}
\begin{lemma}   \label{first}
Let $X$ be a noetherian separated scheme.
Let $Y_0 \supset Y_1 \supset \ldots \supset Y_n$ be
a flag of closed subschemes such that
 $Y_0 = X$ and
$Y_i$
is an ample Cartier divisor on the scheme $Y_{i-1}$
($ 1 \le i \le n$).
Let $J_i$ be the ideal sheaves  on $X$
defining the corresponding subschemes $Y_i$ in $X$.
Let $\sigma = (\eta_0, \ldots, \eta_k) \in S_k$.
Then for any
 $i \le \eta_0$,
for any quasicoherent sheaf
 $\f$ on $X$ we have
\begin{equation}  \label{zamech}
A_{\sigma}  (\f) =   \mathop{\pl}\limits_m
A_{\sigma}(\f / J_i^m \f)  \mbox{.}
\end{equation}
\end{lemma}

\begin{nt}
{  \em
We consider  the sheaf $\f / J_i^m \f$
in~(\ref{zamech})
as the sheaf on the scheme $X$.
The corresponding  functor  of direct image
from the subscheme $Y_i$
is omitted for the sake of simplication of notations.
Further we shall do the same in analogous situations.}
\end{nt}

\noindent
\proof
From the definition of the functor
 $A_{\sigma}$
we have
$$
%\begin{array}{c}
\mathop{\pl}\limits_m
A_{\sigma} (\f / J_i^m \f)  =
\mathop{\pl}\limits_m
\mathop{\pl}\limits_l
A_{\partial_0 (\sigma)}
(
(j_{\eta_0})_*
(j_{\eta_0})^*
(\f / (J_i^m + J_{\eta_0}^l) \f)
)
=   $$
$$
=
\mathop{\pl}\limits_{(m,l)}
A_{\partial_0 (\sigma)}
(
(j_{\eta_0})_*
(j_{\eta_0})^*
(\f / (J_i^m + J_{\eta_0}^l) \f)
)   =
\mathop{\pl}\limits_l
A_{\partial_0 (\sigma)}
(
(j_{\eta_0})_*
(j_{\eta_0})^*
(\f / J_{\eta_0}^l \f)
)
= A_{\sigma} (\f)  \mbox{,}
%\end{array}
$$
where  next to the last equality
follows from cofinality of the projective systems
$$
A_{\partial_0 (\sigma)}
(
(j_{\eta_0})_*
(j_{\eta_0})^*
(\f / (J_i^m + J_{\eta_0}^l) \f)
)
\quad \mbox{and}  \quad
A_{\partial_0 (\sigma)}
(
(j_{\eta_0})_*
(j_{\eta_0})^*
(\f / J_{\eta_0}^l \f)
) \mbox{.}
$$
The last follows from cofinality of the systems
$
\f / (J_i^m + J_{\eta_0}^l) \f
$
and
$
\f / J_{\eta_0}^l \f
$.
Besides, in our reasonings we meant that if
 $\sigma \in S_0 $,
then $A_{\partial_0 (\sigma)} = H^0 (X, \cdot)$.
The lemma is proved.

\begin{lemma}     \label{zv-1}
Let $X$ be a Cohen-Macaulay noetherian scheme.
Let $Y_0 \supset Y_1 \supset \ldots \supset Y_n$ be
a flag of closed subschemes
such that
 $Y_0=X$ and
$Y_i$
is a Cartier divisor on the scheme
$Y_{i-1}$ ($1 \le i \le n$).
Let $J_i$ be the ideal sheaves on $X$  defining
the corresponding
 subschemes $Y_i$ in $X$.
Let $j_k$ be an open imbedding
of $Y_k \backslash Y_{k+1}$
into $Y_k$  ($0 \le k < n$).
Then for any  $0 \le k < n$,
for any $m \ge 1$,
for any locally free sheaf
 $\f$ on $X$
the natural map
$$
\f / J_k^m \f  \lto (j_k)_* (j_k)^*  (\f / J_k^m \f)
$$
is an imbedding.
\end{lemma}
\proof
Let us do induction on $m$.

Let $m = 1$.
Then $\f / J_k \f$ is a locally free sheaf on $Y_k$.
Therefore,
applying lemma~\ref{lem9}
to the pair $Y_{k+1} \hookrightarrow Y_k$,
we obtain injectivity in this case,
since  Cartier divisor is locally
generated by one element
which is not divisor of zero in the structure sheaf.

Let $m >1$.
Then since
$Y_i$ are Cartier divisors on
 $Y_{i-1}$,
we have that
 $Y_k$
is a local complete intersection in $X$
and locally defined by a regular
sequence on $X$.
Therefore the sheaf
 $J_k^{m-1}/ J_k^m$
is locally free on $Y_k$ (see~\cite[ ch.~II, th.~8.21A]{Ha}).
From the exact sequence
\begin{equation}  \label{DEDE}
0 \lto
J_k^{m-1} \f / J_k^m  \f
\lto
\f / J_k^m \f
\lto
\f / J_k^{m-1} \f
\lto
0
\end{equation}
and lemma~\ref{lem5}
we obtain exactness of the following sequence
\begin{equation} \label{DEDEDE}
0
\to
(j_k)_* (j_k)^*
(
J_k^{m-1} \f  / J_k^m  \f
)
\to
(j_k)_* (j_k)^*
(
\f / J_k^m \f
)
\to
(j_k)_* (j_k)^*
(
\f / J_k^{m-1} \f
)
\to
0
\mbox{.}
\end{equation}
The sheaf $
J_k^{m-1} \f  / J_k^m  \f
=
\f
\otimes_{\oo_X}
J_k^{m-1} / J_k^m
$
is locally free $Y_k$.
Therefore the map
$$
J_k^{m-1} \f / J_k^m  \f
\lto
(j_k)_* (j_k)^*
(
J_k^{m-1} \f / J_k^m  \f
)
$$
is an  imbedding.
Also the map
$$
\f / J_k^{m-1} \f  \lto
(j_k)_* (j_k)^* (\f / J_k^{m-1} \f)
$$
is an imbedding by induction hypothesis.
From this, (\ref{DEDE}), (\ref{DEDEDE})
and non complicated diagram search
we obtain that the map
$$
\f / J_k^{m} \f  \lto
(j_k)_* (j_k)^* (\f / J_k^{m} \f)
$$
is an imbedding. The lemma is proved.

\begin{lemma}       \label{zv-}
Let $X$ be a projective equidimensional
Cohen-Macaulay  algebraic scheme of
dimension
$n$  over  a field.
Let $Y_0 \supset Y_1 \supset \ldots \supset Y_k$ be a flag
of closed subschemes such that
$Y_0 = X$ and
$Y_i$
is an ample Cartier divisor on the scheme
 $Y_{i-1}$
for any
$i$ ($1 \le i \le k$).
Let $J_i$ be
the ideal sheaves on $X$ defining the corresponding
subschemes $Y_i$ in $X$.
Then for any locally free sheaf
$\f$
on $X$   the natural map
$$
H^0(X, \f) \lto
\mathop{\pl}\limits_{m }
H^0(X, \f / J_k^m \f)
$$
is an imbedding;
and if $k < n$,
then this one is an isomorphism.
\end{lemma}
\proof
The proof will be done by induction on $k$.

Let at first $k=1$.
The sheaf $J_1$
is the  dual of the ample invertible sheaf on $X$.
And from conditions on  $X$
(that is Cohen-Macaulayness, projectiveness , equidimensionality)
 there exist (see~\cite[ ch.~III, th.~7.6]{Ha})
$l > 0$
such that for any
 $m > l$  we have
$H^0 (X, J_1^m \f) = 0$,  and if
$n \ge 2$, then  we have $H^1(X, J_1^m \f) =0$ as well.

Hence and from the exact sequence
$$
0 \lto J_1^m \f \lto  \f \lto \f/ J_1^m \f \lto 0
$$
we obtain that the map
$$
H^0(X, \f)  \lto
H^0(X, \f / J_1^m \f)
$$
is an imbedding for $m > l$, and
$$
H^0(X, \f)  =
H^0(X, \f / J_1^m \f)
$$
for $m > l$ and $n \ge 2$.
And after passage to the projective limit on $m$
we obtain the imbedding
$$
H^0(X, \f) \hookrightarrow
\mathop{\pl}\limits_{m }
H^0(X, \f / J_1^m \f) =
H^0 (X,
\mathop{\pl}\limits_{m }
\f / J_1^m \f
)     \mbox{,}
$$
and  if $n > 1$, then this one is the isomorphism
$$
H^0(X, \f)  \simeq
H^0 (X,
\mathop{\pl}\limits_{m }
\f / J_1^m \f
) \mbox{.}
$$

Now let $k>1$ be arbitrary.
From lemmas conditions it follows that
 $Y_k$  is locally defined by a regular sequence on $X$.
Therefore $Y_k$ is a Cohen-Macaulay scheme,
and $J_k^m / J_k^{m+1}$ are locally free sheaves on
 $Y_k$ (see~\cite[ ch.~II, th.~8.21A]{Ha}).
By induction hypothesis we have
\begin{equation}     \label{DOM}
\mathop{\pl}\limits_{l }
H^0(X, \f / J_{k-1}^l \f) =
H^0(X, \f)  \mbox{.}
\end{equation}
From cofinality of projective systems
$\f / J_k^m \f$
and
$\f / (J_k^m + J_{k-1}^l) \f$
we have that
$$
\mathop{\pl}\limits_{m}
H^0(X, \f / J_{k}^m \f) =
\mathop{\pl}\limits_{(l,m)}
H^0(X, \f / (J_{k}^m + J_{k-1}^l) \f) =
\mathop{\pl}\limits_{l}
\mathop{\pl}\limits_{m}
H^0(X, \f / (J_k^m + J_{k-1}^l) \f) \mbox{.}
$$
Frim this one and equality~(\ref{DOM})
it suffices for the proof of the lemma to show that
for any $l \ge 1$
the map
$$
H^0(X, \f / J_{k-1}^l \f)
\lto
\mathop{\pl}\limits_{m}
H^0(X, \f / (J_k^m + J_{k-1}^l) \f)
$$
is an imbedding,
and if $k < n $ then this map is an isomorphism.

For this one let us consider the exact sequence
$$
0 \lto
\frac{J_k^m \f + J_{k-1}^l \f}{ J_{k-1}^l \f}
\lto
\f / J_{k-1}^l \f
\lto
\f / (J_k^m + J_{k-1}^l) \f  \lto  0 \mbox{.}
$$
It suffices to show that for all sufficiently large
 $m$
\begin{equation} \label{DI}
H^0(X, \frac{J_k^m \f + J_{k-1}^l \f}{ J_{k-1}^l \f}) = 0   \mbox{,}
\end{equation}
and if $k < n $, then
\begin{equation}     \label{DIDI}
H^1(X, \frac{J_k^m \f + J_{k-1}^l \f}{ J_{k-1}^l \f}) = 0   \mbox{.}
\end{equation}

For this one let us do induction on $l$.
If $l=1$,
then (\ref{DI}) and (\ref{DIDI})
follows at once from~\cite[ ch.~III, th.~7.6]{Ha}
and the fact that the sheaf
$J_k/J_{k-1}$  is the dual of the ample invertible sheaf on $Y_{k-1}$.

If $l>1$, then from the identity $ \frac{A + B}{B} = \frac{A}{A \cap B} $
it follows the exact sequence
$$
0
\lto
\frac{J_k^m   \cap J_{k-1}^{l-1}}{J_k^m \cap J_{k-1}^l}
\otimes_{\oo_X} \f
\lto
\frac{J_k^m \f + J_{k-1}^l \f}{J_{k-1}^l \f}
\lto
\frac{J_k^m \f + J_{k-1}^{l-1} \f}{J_{k-1}^{l-1} \f}  \lto 0  \mbox{.}
$$
Restricting ourself to the local situation and
using the fact that
$J_i$  is generated by  a regular sequence  and~\cite[ ch.~II, th.~8.21A]{Ha},
it is not difficult to understand that
$$
\frac{J_k^m   \cap J_{k-1}^{l-1}}{J_k^m \cap J_{k-1}^l}
\otimes_{\oo_X} \f
=
\frac{J_k^{m-l+1} + J_{k-1}}{J_{k-1}}
\cdot
\left(
\frac{J_{k-1}^{l-1}}{J_{k-1}^l}
\otimes_{\oo_X} \f
\right)  \mbox{.}
$$
In addition, the sheaf
$
\frac{J_{k-1}^{l-1}}{J_{k-1}^l}
\otimes_{\oo_X} \f
$
is locally free on $Y_{k-1}$.
The sheaf $J_k/J_{k-1}$
is the dual of the ample invertible sheaf on the Cohen-Macaulay scheme
$Y_{k-1}$.
Therefore
from~\cite[ ch.~III, th.~7.6]{Ha}
we have for sufficiently large
$m$ that:
$$
H^0 (X,
\frac{J_k^m   \cap J_{k-1}^{l-1}}{J_k^m \cap J_{k-1}^l}
\otimes_{\oo_X} \f ) = 0     \qquad  \mbox{, and}
$$
$$
\mbox{if}  \quad
k < n \quad \mbox{, then}  \qquad
H^1 (X,
\frac{J_k^m   \cap J_{k-1}^{l-1}}{J_k^m \cap J_{k-1}^l}
\otimes_{\oo_X} \f )  = 0     \mbox{.}    \qquad \qquad \qquad \qquad \quad \;
$$
The lemma is proved.

\medskip
\bigskip
\noindent
{\bf Corollary}(from lemma~\ref{zv-}) \\
{\em
Under the conditions of lemma~\ref{zv-}
for any $\sigma \in S_0$
the natural map
$$
H^0(X, \f)  \lto A_{\sigma}(\f)
$$
is an imbedding.  }

\bigskip
\noindent
\proof
Let $\sigma = (m)$.
By lemma~\ref{zv-}
we have the imbedding
$$
0 \lto
\f / J_k^m \f
\lto
(j_k)_*
(j_k)^*
(\f /J_k^m \f)  \mbox{.}
$$
Hence we obtain the imbedding
$$
0  \lto
H^0(X, \f / J_k^m \f)
\lto
H^0 (X,
(j_k)_*
(j_k)^*
(\f /J_k^m \f)
 )   \mbox{.}
$$
After passage to the projective limit on
$m$
we obtain the imbedding
$$
0 \lto
\mathop{\pl}\limits_m
H^0 (X, \f /J_k^m \f)
\lto
A_{\sigma} (\f)  \mbox{.}
$$
Now application of lemma~\ref{zv-}
concludes the proof of the corollary.

\medskip
\begin{th}      \label{tend1}
Let $X$
be a projective equidimensional
Cohen-Macaluay scheme of dimension
$n$  over a field.
Let $Y_0 \supset Y_1 \supset \ldots
\supset Y_n$ be a flag of closed subschemes
such that
 $Y_0 = X$ and
$Y_i$         is an ample Cartier divisor on the scheme
$Y_{i-1}$
for any
 $1 \le i \le n$.
Then for any locally free sheaf
$\f$
on $X$ we have that
\begin{enumerate}
\item     \label{unkt1}
for any $\sigma  \in S_k$  ($0 \le k \le n$)
the natural map
$$
H^0(X, \f) \lto A_{\sigma} (\f)
$$
is an imbedding,
\item    \label{unkt2}
for any $\sigma \in S_k$  ($1 \le k \le n$),
for any $i$ ($0 \le i \le k$)
the natural map
$$
d_i (\sigma)  \: : \:
A_{\partial_i (\sigma)} (\f)  \lto A_{\sigma} (\f)
$$
ia an imbedding.
\end{enumerate}
\end{th}

\begin{nt}
{\em
Taking into account~(\ref{npp})
from proposition~\ref{predl3},
it is possible to reformulate
item~\ref{unkt2} of this theorem
in the following way: \\
{\em
for any locally free sheaf $\f$
on $X$, for any $\sigma_1, \sigma_2 \in S$, $ \sigma_1 \subset \sigma_2$
the natural map
$
\quad  A_{\sigma_1}(\f)  \lto A_{\sigma_2}(\f)  \quad
$
ia an imbedding.
}}
\end{nt}
%\bigskip \\
\proof
Let $J_i$ be the ideal sheaves
on $X$ defining
the corresponding  subschemes $Y_i$ in $X$.
Let us prove first item~\ref{unkt2} of the theorem.
Consider 3 cases.
\medskip             \\
\underline{Case 1.}
Let $\sigma = (\eta_0, \eta_1, \ldots, \eta_k)$,
and $i=0$.
Then $\partial_0 (\sigma)=
(\eta_1, \ldots, \eta_k)$.

By lemma~\ref{zv-1}
for any $m \ge 1$
the map
$$
\f / J_{\eta_0}^m \f
\lto
(j_{\eta_0})_* (j_{\eta_0})^*  (\f / J_{\eta_0}^m \f)
$$
is an imbedding.
Apply to this sequence the exact functor
 $A_{\partial_0 (\sigma)}$.
We obtain the imbedding:
$$
A_{\partial_0 (\sigma)}
(\f / J_{\eta_0}^m \f)
\lto
A_{\partial_0 (\sigma)}
((j_{\eta_0})_*  (j_{\eta_0})^*
(\f / J_{\eta_0}^m \f) ) \mbox{.}
$$
After passage to the projective limit on
$m$
we obtain the imbedding
$$
\mathop{\pl}\limits_m
A_{(\eta_1, \ldots , \eta_k)}
(\f / J_{\eta_0}^m \f)
\lto A_{\sigma} (\f) \mbox{.}
$$
In addition, from lemma~\ref{first}
we have that
$$
\mathop{\pl}\limits_m
A_{(\eta_1, \ldots , \eta_k)}
(\f / J_{\eta_0}^m \f)
=
A_{(\eta_1, \ldots, \eta_k)} (\f)  \mbox{.}
$$
Thus we obtain that in this case the map
$$
A_{\partial_0 (\sigma)} (\f)
\lto
A_{\sigma} (\f)
$$
is an imbedding.
\medskip        \\
\underline{Case 2.}
Let $\sigma = (\eta_0, \eta_1)$,
and $i=1$.
In this case we have that
$$
A_{\partial_1(\sigma)}(\f) =
\mathop{\pl}\limits_m
H^0 (X, (j_{\eta_0})_* (j_{\eta_0})^*  (\f / J_{\eta_0}^m \f ))
$$
$$
A_{\sigma} (\f) =
\mathop{\pl}\limits_m
A_{\partial_0 (\sigma)}
(
(j_{\eta_0})_*  (j_{\eta_0})^*
(\f / J_{\eta_0}^m \f)
)      \mbox{,}
$$
and for the proof of this case  it suffices to show
that for any
$m \ge 1$
the map
\begin{equation} \label{pauk}
H^0 (X, (j_{\eta_0})_*  (j_{\eta_0})^*  (\f / J_{\eta_0}^m \f ) )
\lto
A_{\partial_0 (\sigma)}
((j_{\eta_0})_*    (j_{\eta_0})^*  (\f / J_{\eta_0}^m \f) )
\end{equation}
is an imbedding.

Let us show this by induction.
Let $m=1$.
Then $\f / J_{\eta_0} \f$ is a locally free sheaf on $Y_{\eta_0}$.
Besides,
since $Y_{\eta_0 + 1}$ is a Cartier divisor on
$Y_{\eta_0}$,
we have from item~\ref{pun5} of proposition~\ref{predl1}  that
$$
%\begin{array}{rcl}
(j_{\eta_0})_* (j_{\eta_0})^*
(\f / J_{\eta_0} \f)
 =
((j_{\eta_0})_* (j_{\eta_0})^*
\oo_{Y_{\eta_0}})
\otimes
_{\oo_{Y_{\eta_0}}}
(\f / J_{\eta_0} \f) =   $$
\begin{equation}  \label{ve}
=
\mathop{\il}\limits_j
\oo_{Y_{\eta_0}} (j Y_{\eta_0 +1})
\otimes_{\oo_{Y_{\eta_0}}}
(\f / J_{\eta_0} \f)
 =
\mathop{\il}\limits_j
(\f / J_{\eta_0} \f) (j Y_{\eta_0 +1}) \mbox{.}
%\end{array}
\end{equation}
The sheaves
$(\f / J_{\eta_0} \f) (j Y_{\eta_0 +1}) $
are locally free on
 $Y_{\eta_0}$ as well.

Therefore by corollary from lemma~\ref{zv-} the map
$$
H^0 (X,
(\f / J_{\eta_0} \f) (j Y_{\eta_0 +1})
)
\lto
A_{\partial_0 (\sigma)}
(
(\f / J_{\eta_0} \f) (j Y_{\eta_0 +1})
)
$$
is an imbedding.
After passage to the projective limit on
$j$
we obtain injectivity of~(\ref{pauk})   in the case  $m =1$.

If $m > 1$,
then the statement will follow
from induction hypothesis
and consideration of the following
two exact sequences:
$$
\begin{array}{l}
0 \lto
H^0(X, (j_{\eta_0})_* (j_{\eta_0})^*  \g)
  \lto
H^0(X, (j_{\eta_0})_* (j_{\eta_0})^*  \f / J_{\eta_0}^m
\f)  \lto      \qquad \qquad \qquad     \qquad \qquad \qquad
\end{array}
$$
$$
\begin{array}{r}
\qquad \qquad \qquad \qquad \qquad \qquad \qquad     \qquad \qquad
\qquad
\lto
H^0(X, (j_{\eta_0})_* (j_{\eta_0})^*  \f / J_{\eta_0}^{m-1}
\f)
\lto
0
\end{array}
$$
\smallskip
$$
\begin{array}{l}
0 \lto
A_{\partial_0 (\sigma)}
((j_{\eta_0})_* (j_{\eta_0})^*  \g)
 \lto
A_{\partial_0 (\sigma)}((j_{\eta_0})_* (j_{\eta_0})^*  \f / J_{\eta_0}^m
\f)
\lto          \qquad \qquad \qquad   \qquad \qquad \qquad
\end{array}
$$
$$
\begin{array}{r}
\qquad \qquad \qquad \qquad \qquad     \qquad \qquad \qquad \qquad
\qquad
\lto
A_{\partial_0 (\sigma)}((j_{\eta_0})_* (j_{\eta_0})^*  \f / J_{\eta_0}^{m-1}
\f)
\lto 0 \mbox{,}
\end{array}
$$
where the sheaf $\g = J_{\eta_0}^{m-1} \f / J_{\eta_0}^m \f$
is a locally free on $Y_{\eta_0}$.
 Case~2 is analyzed.
\medskip        \\
\underline{Case 3.}
We shall consider all that are not in cases~1 and~2.
Let $\sigma = (\eta_0, \eta_1, \ldots, \eta_k)$
and $i \ne 0$.
Let us do induction on $i$.
The case $i=0$
is already  analyzed (case~1).
We have
$$
A_{\partial_i (\sigma)} (\f)  =
\mathop{\pl}_m
A_{\partial_{i-1} \cdot (\partial_0 (\sigma))}
(
%i_{\eta_0})_*
(j_{\eta_0})_*  (j_{\eta_0})^* (\f / J_{\eta_0}^m \f)
)      \quad    \mbox{and}
$$
$$
A_{\sigma} (\f)  =
\mathop{\pl}_m
A_{\partial_0 (\sigma)}
(
%i_{\eta_0})_*
(j_{\eta_0})_*  (j_{\eta_0})^* (\f / J_{\eta_0}^m \f)
)   \mbox{.}
$$
By induction hypothesis
applied to the scheme
 $Y_{\eta_0}$
we can suppoce that
for any locally free sheaf
 $\cal H$ on $Y_{\eta_0}$
the map
\begin{equation}     \label{anti}
0 \lto A_{\partial_{i-1} (\partial_0 (\sigma))}  ({\cal H})
\lto
A_{\partial_0 (\sigma)} ({\cal H})
\end{equation}
is an imbedding.
Let us show that for any
 $m$
the map
\begin{equation}    \label{daaaa}
A_{\partial_{i-1} \cdot (\partial_0 (\sigma))}
(
%i_{\eta_0})_*
(j_{\eta_0})_*  (j_{\eta_0})^* (\f / J_{\eta_0}^m \f)
)
\lto
A_{\partial_0 (\sigma)}
(
%i_{\eta_0})_*
(j_{\eta_0})_*  (j_{\eta_0})^* (\f / J_{\eta_0}^m \f)
)
\end{equation}
is an imbedding.

If $m=1$,
then, using~(\ref{ve}),
as in case ~2,
we  reduce all at once   to the sequence~(\ref{anti}).
Afterwards we pass to the direct limit.

If $m>1$, then as in case~2
we can do induction on $m$
by means of using of two following exact sequences:
$$
\begin{array}{l}
0 \lto
A_{\partial_{i-1} (\partial_0 (\sigma))}
(j_{\eta_0})_* (j_{\eta_0})^*  \g)
  \lto
A_{\partial_{i-1} (\partial_0 (\sigma))}
((j_{\eta_0})_* (j_{\eta_0})^*  \f / J_{\eta_0}^m
\f)
\lto   \qquad \qquad \qquad  \qquad \qquad
\end{array}
$$
$$
\begin{array}{r}
\quad
\qquad \qquad \qquad	 \qquad \qquad \qquad  \qquad \qquad  \lto
A_{\partial_{i-1} (\partial_0 (\sigma))}
((j_{\eta_0})_* (j_{\eta_0})^*  \f / J_{\eta_0}^{m-1}
\f)
\lto
0
\end{array}
$$
\smallskip
$$
\begin{array}{l}
 0 \lto
A_{\partial_0 (\sigma)}
((j_{\eta_0})_* (j_{\eta_0})^*  \g)
 \lto
A_{\partial_0 (\sigma)}((j_{\eta_0})_* (j_{\eta_0})^*  \f / J_{\eta_0}^m
\f)
\lto  \qquad \qquad \qquad \qquad \qquad
\end{array}
$$
$$
\begin{array}{r}
\qquad \qquad
\qquad \qquad \qquad \qquad \qquad \qquad  \qquad \qquad \lto
A_{\partial_0 (\sigma)}((j_{\eta_0})_* (j_{\eta_0})^*  \f / J_{\eta_0}^{m-1}
\f)
\lto 0 \mbox{,}
\end{array}
$$
where
the sheaf
  $\g = J_{\eta_0}^{m-1} \f / J_{\eta_0}^m \f$.

Now after passage in~(\ref{daaaa})  to the   limit
we conclude the proof of case~3.

\bigskip
Now consider item~1 of the theorem.
If $k=0$,
then
this is corollary  of lemma~\ref{zv-}.
If $k > 0$,
then consider the map
$$
H^0(X, \f)  \lto A_{\sigma} (\f)
$$
as composition of the maps
$$
H^0(X, \f)   \lto A_{\partial_0 (\sigma)} (\f)      \qquad \mbox{and}
\qquad
A_{\partial_0 (\sigma)}(\f)
\lto
A_{\sigma} (\f)   \mbox{,}
$$
where we can suppose that by induction on $k$ the first map
 is injective, and  by item~2 of this theorem the second
map is injective as well.
Theorem~\ref{tend1}  is proved.

\begin{lemma} \label{de}
Let all the conditions of  theorem~\ref{tend1}
be satisfied.
By
$J_i$
denote the ideal sheaves on $X$ defining
the corresponding $Y_i$ in $X$.
Then for any locally free sheaf $\f$
on $X$, for any $m > 0$ we have that
\begin{enumerate}
\item                \label{pa}
the map
$$
H^0 (X, \f / J_i^m \f)  \lto
A_{\sigma} (\f / J_i^m \f)
$$
is an imbedding for any $\sigma = (\zeta_0, \ldots, \zeta_k)$,
$0 \le i \le \zeta_0$.
\item          \label{pb}
the map
$$
A_{\sigma_1} (\f / J_i^m  \f)
\lto
A_{\sigma_2} (\f / J_i^m \f)
$$
is an imbedding
for any $\sigma_1, \sigma_2 \in S$,
$ \sigma_1 \subset \sigma_2 = (\eta_0, \ldots, \eta_k)$ ,
$0 \le i \le \eta_0$.
\end{enumerate}
\end{lemma}
\proof
Let us show item~\ref{pa}.
If $m =1$,
then this follows from theorem~\ref{tend1},
which is applied to the scheme $Y_i$.

If $m >1$,
then apply induction.
For this consider the exact sequence
$$
0
\lto
J_i^{m-1}  \f  / J_i^m \f
\lto
\f / J_i^m \f
\lto
\f/ J_i^{m-1} \f
\lto
0   \mbox{.}
$$
Hence and from exactness of the functor
$A_{\sigma}$
we have the following commutative diagram:
{
$$
\begin{array}{ccccccc}
0 \lto &
H^0 (X, J_i^{m-1}  \f  / J_i^m \f)
 & \lto
&
H^0 (X, \f / J_i^m \f)
 &
\lto &
H^0 (X, \f/ J_i^{m-1} \f)
  &
\\           &
 \begin{picture}(0,28)
 \put(0,26){\vector(0,-1){26}}
 \put(0,11){$\; \alpha $}
 \end{picture}
 &&
 \begin{picture}(0,28)
 \put(0,26){\vector(0,-1){26}}
 \put(0,10){$\; \beta $}
 \end{picture}
 &&
 \begin{picture}(0,28)
 \put(0,26){\vector(0,-1){26}}
 \put(0,11){$\; \gamma $}
 \end{picture}
         \\
 0 \lto
 & A_{\sigma} (J_i^{m-1}  \f  / J_i^m \f)
     &
\lto &
A_{\sigma} (\f / J_i^m \f)
&
\lto
&
A_{\sigma}(\f/ J_i^{m-1} \f)
&
\lto 0 \mbox{.}
\end{array}
$$
}
The sheaf $J_i^{m-1}  \f  / J_i^m \f=
\frac{J_i^{m-1}}{J_i^m} \otimes_{\oo_X} \f$ is
locally free on $Y_i$ (see~\cite[ ch.~II, th.~8.21A]{Ha}).
Therefore the map
 $\alpha$ is an imbedding.
The map $\gamma$ is an imbedding
by the induction hypothesis.
Now from non complicated diagram search
it follows that $\beta$
is an imbedding as well. Item~\ref{pa}
is proved.

Item~\ref{pb}  follows at once
from analogous to item~\ref{pa}
and inductive on $m$ reasonings
applied to the following diagram:
{
$$
\begin{array}{ccccccc}
0 \lto &
A_{\sigma_1} (J_i^{m-1}  \f  / J_i^m \f)
 & \lto
&
A_{\sigma_1} ( \f / J_i^m \f)
 &
\lto &
A_{\sigma_1} (\f/ J_i^{m-1} \f)
  &
\\           &
 \begin{picture}(0,28)
 \put(0,26){\vector(0,-1){26}}
 \put(0,11){$\; \alpha $}
 \end{picture}
 &&
 \begin{picture}(0,28)
 \put(0,26){\vector(0,-1){26}}
 \put(0,10){$\; \beta $}
 \end{picture}
 &&
 \begin{picture}(0,28)
 \put(0,26){\vector(0,-1){26}}
 \put(0,11){$\; \gamma $}
 \end{picture}
         \\
 0 \lto
 & A_{\sigma_2} (J_i^{m-1}  \f  / J_i^m \f)
     &
\lto &
A_{\sigma_2} (\f / J_i^m \f)
&
\lto
&
A_{\sigma_2}(\f/ J_i^{m-1} \f)
&
\lto 0 \mbox{.}
\end{array}
$$
}
The lemma is proved.

\bigskip
\begin{th}       \label{tend2}
Let all the conditions of theorem~\ref{tend1}
be satisfied.
Then for any locally free sheaf $\f$,
for any $\sigma_1, \sigma_2 \in S$ we have that
\begin{enumerate}
\item                \label{pu1}
if $\sigma_1 \cap \sigma_2 = \o$,
then
$$
A_{\sigma_1} (\f)  \cap A_{\sigma_2} (\f) = H^0 (X, \f)  \mbox{;}
$$
\item                 \label{pu2}
if $\sigma_1 \cap \sigma_2 \ne \o$,
then
$$
A_{\sigma_1}(\f) \cap A_{\sigma_2} (\f)
= A_{\sigma_1 \cap \sigma_2} (\f)  \mbox{.}
$$
\end{enumerate}
\end{th}
\begin{nt}  {\em
According to theorem~\ref{tend1},
the intersections make sense,
because we can always imbed
$A_{\sigma_1} (\f)$
and $A_{\sigma_2} (\f)$
into $A_{\eta} (\f)$,
where $\eta$
contains $\sigma_1$  and $\sigma_2$.
For instance, $\eta = \sigma_1 \cup \sigma_2$.}
\end{nt}
\noindent
\proof
As usually,
by
 $J_i$  denote  the ideal sheaves on $X$  defining
the corresponding subschemes $Y_i$ in $X$.

Let us show item~\ref{pu1}.
Let $\sigma_1 = (\eta_0, \ldots)$,
$\sigma_2 = (\zeta_0, \ldots)$,
$\sigma_1 \cap \sigma_2 = \o$.
Without loss of generality it can be assumed that
 $\zeta_0 > \eta_0$.
Assume that $\sigma_1 \notin S_0$.

By lemma~\ref{zv-1}
for any $m >1$
we have the exact sequence:
$$
0
\lto
\f / J_{\eta_0}^m \f
\lto
(j_{\eta_0})_* (j_{\eta_0})^*
(\f / J_{\eta_0}^m \f)
\lto
\g_m
\lto
0   \mbox{,}
$$
where the sheaf
$
\g_m  = \frac{
(j_{\eta_0})_* (j_{\eta_0})^*
(\f / J_{\eta_0}^m \f)
}{
\f / J_{\eta_0}^m \f
}
$.

Let us show by induction on $m$
that the natural map
\begin{equation}      \label{DAR}
A_{\partial_0 (\sigma_1)}(\g_m)
\lto
A_{\partial_0 (\sigma_1) \cup \sigma_2} (\g_m)
\end{equation}
is an imbedding.

If $m=1$, then
$$
\g_1 =
\frac{
(j_{\eta_0})_* (j_{\eta_0})^*
(\f / J_{\eta_0} \f)
}{
\f / J_{\eta_0} \f
}    =
\frac{
(j_{\eta_0})_* (j_{\eta_0})^*
\oo_{Y_{\eta_0}}
\otimes_{\oo_{Y_{\eta_0}}}
(\f / J_{\eta_0} \f)
}{\f / J_{\eta_0} \f} =
$$
$$
=
\frac{ \mathop{\il}\limits_k
\oo (kY_{\eta_0 +1})
\otimes_{\oo_{Y_{\eta_0}}}    (\f / J_{\eta_0} \f)
}{\f / J_{\eta_0} \f}=
\mathop{\il}\limits_{k > 0}
(
(\f / J_{\eta_0} \f)  \otimes_{\oo_{Y_{\eta_0}}}
(\oo_{Y_{\eta_0}}(k Y_{\eta_0 +1}) /
\oo_{Y_{\eta_0}}  )
)  \mbox{.}
$$
Denote the sheaf $\h_k=
(\f / J_{\eta_0} \f)  \otimes_{\oo_{Y_{\eta_0}}}
(\oo_{Y_{\eta_0}}(k Y_{\eta_0 +1}) /
\oo_{Y_{\eta_0}} )
$.
By induction on $k$
let us show that the maps
\begin{equation}  \label{DARA}
A_{\partial_0 (\sigma_1)} (\h_k)
\lto
A_{\partial_0 (\sigma_1) \cup \sigma_2} (\h_k)
\end{equation}
are imbeddings.

If $k=1$,
then the sheaf $\h_1$
is locally free on $Y_{\eta_0 +1}$.
Therefore in this case~(\ref{DARA})
follows from theorem~\ref{tend1}
applied to $Y_{\eta_0 +1}$.

If $k >1$,
then from the  exact sequence
$$
0 \lto
\frac{ \oo_{Y_{\eta_0}} ( (k-1) Y_{\eta_0 +1} ) }{
\oo_{Y_{\eta_0}}
}
\lto
\frac{\oo_{Y_{\eta_0}} ( k Y_{\eta_0 +1} ) }{
\oo_{Y_{\eta_0}} }
\lto
\frac{\oo_{Y_{\eta_0}}( k Y_{\eta_0 +1}) }{
\oo_{Y_{\eta_0}}( (k-1) Y_{\eta_0 +1})
}
\lto
0
$$
it follows the commutative diagram:
{
$$
\begin{array}{ccccccc}
0 \lto &
A_{\partial_0 (\sigma_1)} (\h_{k-1})
 & \lto
&
A_{\partial_0 (\sigma_1)} (\h_{k})
 &
\lto &
A_{\partial_0(\sigma_1)} (
\h_k / \h_{k-1}
)
&
\lto
0

\\           &
 \begin{picture}(0,28)
 \put(0,26){\vector(0,-1){26}}
 \put(0,11){$\; \alpha $}
 \end{picture}
 &&
 \begin{picture}(0,28)
 \put(0,26){\vector(0,-1){26}}
 \put(0,10){$\; \beta $}
 \end{picture}
 &&
 \begin{picture}(0,28)
 \put(0,26){\vector(0,-1){26}}
 \put(0,11){$\; \gamma $}
 \end{picture}
         \\
 0 \lto
 &
A_{\partial_0 (\sigma_1) \cup \sigma_2} (\h_{k-1})
 &
\lto &
A_{\partial_0 (\sigma_1) \cup \sigma_2} (\h_k)
&
\lto
&
A_{\partial_0 (\sigma_1) \cup \sigma_2} (\h_k / \h_{k-1})
&
\lto 0 \mbox{.}
\end{array}
$$
}
The sheaf $\h_k / \h_{k-1} =
(\f / J_{\eta_0} \f)  \otimes_{\oo_{Y_{\eta_0}}}
(\oo_{Y_{\eta_0}}(k Y_{\eta_0 +1}) /
\oo_{Y_{\eta_0}}((k-1) Y_{\eta_0 +1}) )
$
is locally free on $Y_{\eta_0 +1}$.
Therefore the map
 $\gamma$
is injective by theorem~\ref{tend1}
applied to $Y_{\eta_0 +1}$.
The map $\alpha$
is injective by induction hypothesis.
Hence we have that the map $\beta$ is injective as well.
Thus (\ref{DARA}) is proved.
After passage in~(\ref{DARA})
to the direct limit on $k$
we obtain~(\ref{DAR}) in the case $m=1$.

Now let us show~(\ref{DAR})
in the case $m > 1$.
From the  exact sequences:
$$
0 \lto
\frac{J_{\eta_0}^{m-1} \f}{ J_{\eta_0}^m \f}
\lto
\f / J_{\eta_0}^m \f
\lto
\f / J_{\eta_0}^{m-1} \f
\lto
0
$$
and
$$
0 \lto
(j_{\eta_0})_* (j_{\eta_0})^*
(\frac{J_{\eta_0}^{m-1} \f}{ J_{\eta_0}^m \f})
\lto
(j_{\eta_0})_* (j_{\eta_0})^*
(\f / J_{\eta_0}^m \f)
\lto
(j_{\eta_0})_* (j_{\eta_0})^*
(\f / J_{\eta_0}^{m-1} \f)
\lto
0
$$
it follows the exact sequence
$$
0
\lto
\frac{
(j_{\eta_0})_* (j_{\eta_0})^*
(\frac{J_{\eta_0}^{m-1} \f}{ J_{\eta_0}^m \f})
}{
\frac{J_{\eta_0}^{m-1} \f}{ J_{\eta_0}^m \f}
}
\lto
\g_m
\lto
\g_{m-1}
\lto
0  \mbox{.}
$$
Applying the exact functors
$A_{\partial_0 (\sigma_1)}$
and
$A_{\partial_0 (\sigma_1) \cup \sigma_2}$,
we obtain the  diagram
{
$$
\begin{array}{ccccccc}
0 \lto &
A_{\partial_0 (\sigma_1)} (
\frac{
(j_{\eta_0})_* (j_{\eta_0})^*
(\frac{J_{\eta_0}^{m-1} \f}{ J_{\eta_0}^m \f})
}{
\frac{J_{\eta_0}^{m-1} \f}{ J_{\eta_0}^m \f}
}
)
 & \lto
&
A_{\partial_0 (\sigma_1)} (\g_m)
 &
\lto &
A_{\partial_0(\sigma_1)} (\g_{m-1})
\lto
0
  &
\\           &
 \begin{picture}(0,36)
 \put(0,32){\vector(0,-1){32}}
 \put(0,14){$\; \alpha $}
 \end{picture}
 &&
 \begin{picture}(0,36)
 \put(0,32){\vector(0,-1){32}}
 \put(0,14){$\; \beta $}
 \end{picture}
 &&
 \begin{picture}(0,36)
 \put(0,32){\vector(0,-1){32}}
 \put(0,14){$\; \gamma $}
 \end{picture}
         \\
 0 \lto
 &
A_{\partial_0 (\sigma_1) \cup \sigma_2} (
\frac{
(j_{\eta_0})_* (j_{\eta_0})^*
(\frac{J_{\eta_0}^{m-1} \f}{ J_{\eta_0}^m \f})
}{
\frac{J_{\eta_0}^{m-1} \f}{ J_{\eta_0}^m \f}
}
)
 &
\lto &
A_{\partial_0 (\sigma_1) \cup \sigma_2} (\g_m)
&
\lto
&
A_{\partial_0 (\sigma_1) \cup \sigma_2} (\g_{m-1})
\lto 0 \mbox{.}
\end{array}
$$
}
The sheaf $\frac{J_{\eta_0}^{m-1} \f}{ J_{\eta_0}^m \f} =
\frac{J_{\eta_0}^{m-1}}{J_{\eta_0}^m}
\otimes_{Y_{\eta_0}}  (\f / J_{\eta_0} \f)
$
is locally free on $Y_{\eta_0}$.
Therefore the map $\alpha$
is injective by the same reasons as in the case
 $m=1$.
The map $\gamma$
is injective by the inductive hypothesis.
Therefore the map $\beta$
is injective as well (this follows
from non complicated diagram search).
Thus we have shown~(\ref{DAR}).

Now consider the following diagram.
{
$$
\begin{array}{ccccc}
0 \to
A_{\partial_0 (\sigma_1)} (
\f / J_{\eta_0}^m \f
)
 & \to
&
A_{\partial_0 (\sigma_1)}
(
(j_{\eta_0})_* (j_{\eta_0})^*   (\f / J_{\eta_0}^m \f)
)
 &
\stackrel{\psi}{\to} &
A_{\partial_0(\sigma_1)} (\g_m)
\to
0
 \\
 \begin{picture}(0,28)
 \put(0,26){\vector(0,-1){26}}
 \put(0,11){$\; \phi $}
 \end{picture}
 &&
 \begin{picture}(0,28)
 \put(0,26){\vector(0,-1){26}}
 \put(0,10){$\; \theta $}
 \end{picture}
 &&
 \begin{picture}(0,28)
 \put(0,26){\vector(0,-1){26}}
 \put(0,11){$\; \beta $}
 \end{picture}
         \\
 0 \to
A_{\partial_0 (\sigma_1) \cup \sigma_2} (
\f / J_{\eta_0}^m \f
)
 &
\to &
A_{\partial_0 (\sigma_1) \cup \sigma_2}
(
(j_{\eta_0})_* (j_{\eta_0})^*   (\f / J_{\eta_0}^m \f)
)
&
\to
&
A_{\partial_0 (\sigma_1) \cup \sigma_2} (\g_{m})
\to 0 \mbox{.}
\end{array}
$$
}
(Note also that here the map $\theta$ is injective.
This follows from the fact that
the map  $\beta$  is injective by  the above,
and the statement $\phi$
is injective by lemma~\ref{de}.)

Now let an element
$$
x \in
A_{\partial_0 (\sigma_1)}
(
(j_{\eta_0})_* (j_{\eta_0})^*   (\f / J_{\eta_0}^m \f)
)
\mbox{,  \quad but}
$$
$$
x \notin
A_{\partial_0 (\sigma_1)} (
\f / J_{\eta_0}^m \f
)  \mbox{.}
$$

Then since the map
 $\beta$ is injective,
we have that
\begin{equation}  \label{kru}
\beta \psi (x) \ne 0  \mbox{.}
\end{equation}
Now consider the diagram
{
$$
\begin{array}{ccccc}
0 \to
A_{\sigma_2} (
\f / J_{\eta_0}^m \f
)
 & \to
&
A_{\sigma_2}
(
(j_{\eta_0})_* (j_{\eta_0})^*   (\f / J_{\eta_0}^m \f)
)
 &
\stackrel{\psi_1}{\to} &
A_{\sigma_2} (\g_m)
\to
0
 \\
 \begin{picture}(0,28)
 \put(0,26){\vector(0,-1){26}}
% \put(0,14){$\; \phi $}
 \end{picture}
 &&
 \begin{picture}(0,28)
 \put(0,26){\vector(0,-1){26}}
% \put(0,14){$\; \theta $}
 \end{picture}
 &&
 \begin{picture}(0,28)
 \put(0,26){\vector(0,-1){26}}
 \put(0,11){$\; \beta_1 $}
 \end{picture}
         \\
 0 \to
A_{\partial_0 (\sigma_1) \cup \sigma_2} (
\f / J_{\eta_0}^m \f
)
 &
\to &
A_{\partial_0 (\sigma_1) \cup \sigma_2}
(
(j_{\eta_0})_* (j_{\eta_0})^*   (\f / J_{\eta_0}^m \f)
)
&
\to
&
A_{\partial_0 (\sigma_1) \cup \sigma_2} (\g_{m})
\to 0 \mbox{.}
\end{array}
$$
}
(Similarly to the previous reasonings
we have that in this diagram all the vertical arrows
are injective.)

And if an element
$$
x \in A_{\sigma_2} (\f / J_{\eta_0}^m \f) \mbox{, \quad then}
$$
\begin{equation} \label{kru1}
\beta_1 \psi_1 (x) = 0     \mbox{.}
\end{equation}
Now if
$$
x \in
A_{\partial_0 (\sigma_1)}
(
(j_{\eta_0})_* (j_{\eta_0})^*   (\f / J_{\eta_0}^m \f)
)
\cap
A_{\sigma_2} (\f / J_{\eta_0}^m \f)
$$
(where the intersection  is possible to be taken in
$
A_{\partial_0 (\sigma_1) \cup \sigma_2}
(
(j_{\eta_0})_* (j_{\eta_0})^*   (\f / J_{\eta_0}^m \f)
)
$, because all the appearing maps are injective),
then from functoriality we have
$$
\beta \psi (x)
=
\beta_1 \psi_1 (x)  \mbox{.}
$$
Therefore, comparing this with~(\ref{kru})
and~(\ref{kru1}),
we obtain that
in  \newline
$
A_{\partial_0 (\sigma_1) \cup \sigma_2}
(
(j_{\eta_0})_* (j_{\eta_0})^*   (\f / J_{\eta_0}^m \f)
)
$
is satisfied
$$
A_{\partial_0 (\sigma_1)}
(
(j_{\eta_0})_* (j_{\eta_0})^*   (\f / J_{\eta_0}^m \f)
)
\cap
A_{\sigma_2} (\f / J_{\eta_0}^m \f)
=
A_{\partial_0 (\sigma_1)} (\f / J_{\eta_0}^m \f)
\cap
A_{\sigma_2} (\f / J_{\eta_0}^m \f)  \mbox{.}
$$
Now from the definition of $A_{\sigma}$ we have that
$$
A_{\sigma_1} (\f) =
\mathop{\pl}_m
A_{\partial_0 (\sigma_1)}
(
(j_{\eta_0})_* (j_{\eta_0})^*   (\f / J_{\eta_0}^m \f)
)           \mbox{,}
$$
from lemma~\ref{first} we have that
$$
A_{\sigma_2} (\f) =
\mathop{\pl}_m A_{\sigma_2} (\f / J_{\eta_0}^m \f)  \mbox{.}
$$
Therefore,
$$
A_{\sigma_1} (\f)  \cap A_{\sigma_2} (\f) =
\left(
\mathop{\pl}_m
A_{\partial_0 (\sigma_1)}
(
(j_{\eta_0})_* (j_{\eta_0})^*   (\f / J_{\eta_0}^m \f)
)
\right) \cap
\left(
\mathop{\pl}_m A_{\sigma_2} (\f / J_{\eta_0}^m \f)
\right) =
$$
$$
=  \mathop{\pl}_m     \left(
A_{\partial_0 (\sigma_1)}
(
(j_{\eta_0})_* (j_{\eta_0})^*   (\f / J_{\eta_0}^m \f)
)
\cap
A_{\sigma_2} (\f / J_{\eta_0}^m \f)
\right) = $$
$$ =
\mathop{\pl}_m
\left(
A_{\partial_0 (\sigma_1)} (\f / J_{\eta_0}^m \f)
\cap
A_{\sigma_2} (\f / J_{\eta_0}^m \f)
\right) =
$$
$$
= \left(
\mathop{\pl}_m
A_{\partial_0 (\sigma_1)} (\f / J_{\eta_0}^m \f)
\right)
\cap
\left(
\mathop{\pl}_m
A_{\sigma_2} (\f / J_{\eta_0}^m \f)
\right)  =
A_{\partial_0 (\sigma_1)} (\f) \cap A_{\sigma_2} (\f)  \mbox{.}
$$
Acting further  in this manner,
i.~e.,   eliminating
the minmal number in the union of indices every time,
we obtain that
$$
A_{\sigma_1}(\f) \cap
A_{\sigma_2}(\f)
=
A_{(i)} (\f)  \cap
A_{\sigma} (\f)  \mbox{, \quad where}
$$
$\sigma= (\zeta_0, \ldots)$
and $\zeta_0 > i$.
But in this case by the reasonings, which is completely  analogous
to the above,
we obtain at once that
$$
H^0 (X, (j_i)_* (j_i)^* (\f / J_i^m \f ))
\cap
A_{\sigma} (\f  /  J_i^m \f) =
$$
$$
=
H^0(X, \f /  J_I^m \f)
\cap A_{\sigma} (\f / J_i^m \f)=
H^0 (X, \f / J_i^m \f)   \mbox{.}
$$
(Note that in contrast to the reasonings above with
the functor
$A_{\gamma}$,
the functor $H^0 (X, \cdot)$
is a left exact functor only.
But the key diagram works in this case as well:
{
$$
\begin{array}{ccccccc}
0 \lto &
A_1
 & \lto
&
A_2
 &
\lto &
A_3

  &
\\           &
 \begin{picture}(0,24)
 \put(0,23){\vector(0,-1){23}}
 \put(0,10){$\; \alpha $}
 \end{picture}
 &&
 \begin{picture}(0,24)
 \put(0,23){\vector(0,-1){23}}
 \put(0,9){$\; \beta $}
 \end{picture}
 &&
 \begin{picture}(0,24)
 \put(0,23){\vector(0,-1){23}}
 \put(0,10){$\; \gamma $}
 \end{picture}
         \\
 0 \lto
 &
B_1
 &
\lto &
B_2
&
\lto
&
B_3   \mbox{.}
\end{array}
$$
}
If the maps
 $\alpha$
and $\gamma$ are injective,
then the map $\beta$
is injective as well.)

Now
$$
A_{i} (\f)  \cap A_{\sigma} (\f) =
\left(\mathop{\pl}_m
H^0 (X, (j_i)_* (j_i)^* (\f / J_i^m \f ))
\right)     \cap
\left(\mathop{\pl}_m
A_{\sigma} (\f  /  J_i^m \f)
\right) =
$$
$$
=
\mathop{\pl}_m
\left(
H^0 (X, (j_i)_* (j_i)^* (\f / J_i^m \f ))
\cap
A_{\sigma} (\f  /  J_i^m \f)
\right)
=
\mathop{\pl}_m
H^0 (X, \f / J_i^m \f) =
H^0 (X, \f)  \mbox{,}
$$
where the last equality follows from lemma~\ref{zv-}.
Item~\ref{pu1}  of theorem~\ref{tend2}
is proved.

\bigskip
Now let us show item~\ref{pu2}  of the theorem.
Consider a few cases.
\medskip        \\
\underline{Case 1.}
$\sigma_1 \cap \sigma_2 \ne \o$,
$0 \notin \sigma_1$,
$0 \notin \sigma_2$.

By lemma~\ref{first} we have that
$$
A_{\sigma_1}(\f)
=
\mathop{\pl}_m A_{\sigma_1} (\f / J_1^m \f)     \mbox{,}
\qquad     \qquad
A_{\sigma_2}(\f)
=
\mathop{\pl}_m A_{\sigma_2} (\f / J_1^m \f)     \mbox{,}
$$
$$
A_{\sigma_1 \cap \sigma_2} (\f)
=
\mathop{\pl}_m A_{\sigma_1 \cap \sigma_2} (\f / J_1^m \f)     \mbox{.}
$$
Let us show that for any
 $m \ge 1$
\begin{equation}  \label{kol}
A_{\sigma_1 \cap \sigma_2} (\f / J_1^m \f) =
A_{\sigma_1} (\f / J_1^m \f)  \cap
A_{\sigma_2} (\f / J_1^m \f) \mbox{,}
\end{equation}
where the last intersection
is regarded in
 $A_{\sigma_1 \cup \sigma_2} (\f / J_1^m \f)$.
(By lemma~\ref{de}
we can imbed these groups  there.)

Let us prove~(\ref{kol})
by induction on $m$.
Let $m=1$.
In this case
 $\f / J_1 \f$ is a locally free sheaf on $Y_1$.
And equality~(\ref{kol})
turns into the analogous equality~(\ref{kol})
on $Y_1$.
The scheme $Y_1$ has lesser dimension than dimension of $X$.
Applying induction on dimension of scheme,
we can suppose that theorem~\ref{tend2}
is already true for schemes of lesser dimension.
(For schemes of dimension~1
theorem~\ref{tend2}  follows at once
from theorem~\ref{tend1} and theorem~\ref{teorem3}
by trivial reasons.)
Therefore~(\ref{kol})  is true when $m=1$.

Let $m >1$.
Then the exact sequence
$$
0
\lto
J_1^{m-1} \f  / J_1^m \f
\lto
\f / J_1^m \f
\lto
\f / J_1^{m-1} \f
\lto
0
$$
induces the following commutative diagram:
{
$$
\begin{array}{ccccc}
0 \to
A_{\sigma_1 \cap \sigma_2}
(
J_1^{m-1} \f  / J_1^m \f
)
 & \to
&
A_{\sigma_1 \cap \sigma_2}
(\f / J_1^m \f)
 &
\to &
A_{\sigma_1 \cap \sigma_2}
(\f / J_1^{m-1} \f)
\to 0
\\
 \begin{picture}(0,36)
 \put(0,32){\vector(0,-1){32}}
 \put(0,14){$\; \alpha $}
 \end{picture}
 &&
 \begin{picture}(0,36)
 \put(0,32){\vector(0,-1){32}}
 \put(0,14){$\; \beta $}
 \end{picture}
 &&
 \begin{picture}(0,36)
 \put(0,32){\vector(0,-1){32}}
 \put(0,14){$\; \gamma $}
 \end{picture}
         \\
 0 \to
A_{\sigma_1}
(
\frac{J_1^{m-1} \f}{ J_1^m \f}
)
\cap
A_{\sigma_2}
(
\frac{J_1^{m-1} \f}{ J_1^m \f}
)
 &
\to &
H_m
&
\to
&
\frac{H_m}{
A_{\sigma_1}
\bigl(
\frac{J_1^{m-1} \f}{ J_1^m \f}
\bigr)
\cap
A_{\sigma_2}
\bigl(
\frac{J_1^{m-1} \f}{ J_1^m \f}
\bigr)
}
\to 0
   \mbox{,}
\end{array}
$$
}
where $H_m = A_{\sigma_1} (\f / J_1^m \f)
\cap
A_{\sigma_2} (\f /J_1^m \f)   \mbox{.}
$
From $\sigma_1 \cap \sigma_2 \ne \o$
it follows that the functor    $A_{\sigma_1 \cap \sigma_2}$
is exact. Therefore
the upper row of the diagram is  exact.

There is the natural map $\theta$:
$$
\frac{H_m}{
A_{\sigma_1}
\bigl(
\frac{J_1^{m-1} \f}{ J_1^m \f}
\bigr)
\cap
A_{\sigma_2}
\bigl(
\frac{J_1^{m-1} \f}{ J_1^m \f}
\bigr)
}
\lto
H_{m-1}  \mbox{.}
$$
And from the  exact sequences
$$
0
\lto
A_{\sigma_1} (J_1^{m-1} \f  / J_1^m \f)
\lto
A_{\sigma_1} (\f / J_1^m \f)
\lto
A_{\sigma_1 } (\f / J_1^{m-1} \f )
$$
$$
0 \lto
A_{\sigma_2} (J_1^{m-1} \f  / J_1^m \f)
\lto
A_{\sigma_2} (\f / J_1^m \f)
\lto
A_{\sigma_2} (\f / J_1^{m-1} \f )
$$
it follows at once that the map
 $\theta$  is an imbedding.
Besides,
the map $\theta \cdot \gamma$
is the natural map from
$A_{\sigma_1 \cap \sigma_2}  (\f / J_1^{m-1} \f)$
to $H_{m-1}$;
and, consequently,
by the induction hypothesis it is possible to suppose
that $\theta \cdot \gamma$ is an isomorphism.

From the last two facts we obtain at once
that $\gamma$ is an isomorphism.
Since the sheaf $J_1^{m-1} \f / J_1^m \f =
J_1^{m-1} / J_1^m  \otimes_{\oo_X} \f
$
is locally free on $Y_1$
and $\dim Y_1 < \dim X$,

we have  that the map $\alpha$ is
an isomorphism as well.
Therefore from this commutative diagram
it follows that the map $\beta$ is an isomorphism as well.

Thus, equality~(\ref{kol})
is proved.
Now passage in~(\ref{kol})
to the projective limit on $m$
concludes the proof of case~1.
\medskip        \\
\underline{Case 2.}
$0 \in \sigma_1$,
$0 \notin \sigma_2$ (or vice versa),
$\sigma_1 \cap \sigma_2 \ne \o$.

Now by the analogous reasonings,
as in the proof of item~\ref{pu1} of this theorem,
we  obtain at once the following
$$
A_{\sigma_1} (\f) \cap
A_{\sigma_2} (\f)  =
A_{\partial_0 (\sigma_1)} (\f)
\cap
A_{\sigma_2} (\f)  \mbox{;}
$$
that reduces this case
to the case~1, analyzed above.
\medskip        \\
\underline{Case 3.}
$0 \in \sigma_1$,
$0 \in \sigma_2 $.

Then
$$
A_{\sigma_1}(\f)=
\mathop{\il}_{k}
A_{\partial_0 (\sigma_1)} (\f (kY_1)),
\quad
\quad
A_{\sigma_2}(\f)=
\mathop{\il}_{k}
A_{\partial_0 (\sigma_2)} (\f (kY_1))
$$
Now from case~1 (or if  $
\partial_0 (\sigma_1) \cap \partial_0 (\sigma_2) = \o
 $, then from item~\ref{pu1} of this theorem)  we have that
$$
A_{\partial_0 (\sigma_1)} (\f (k Y_1))
\cap
A_{\partial_0 (\sigma_2)}  (\f (k Y_1)) =
A_{\partial_0 (\sigma_1) \cap \partial_0 (\sigma_2)}
(\f (kY_1)) \mbox{.}
$$
(Here $A_{\o}(\cdot) = H^0 (X,\cdot)$).

Therefore,
$$
A_{\sigma_1}(\f)
\cap
A_{\sigma_2} (\f)
=
\left(
\mathop{\il}_{k}
A_{\partial_0 (\sigma_1)} (\f (kY_1))
\right)
\cap
\left(
\mathop{\il}_{k}
A_{\partial_0 (\sigma_2)} (\f (kY_1))
\right)  =
$$
$$
=
\mathop{\il}_k
\left(
A_{\partial_0 (\sigma_1)} (\f (k Y_1))
\cap
A_{\partial_0 (\sigma_2)}  (\f (k Y_1))
\right) =
\mathop{\il}_k
A_{\partial_0 (\sigma_1) \cap \partial_0 (\sigma_2)}
(\f (kY_1)) =
A_{\sigma_1 \cap \sigma_2} (\f) \mbox{.}
$$
Theorem~\ref{tend2} is proved.

\bigskip
\bigskip
In the sequel we shall assume
that all the conditions of theorem~\ref{tend1}
are satisfied,
and a field $k$ is the field of definition of the scheme $X$.
Also, let us assume that
$Y_n=x$, where $x$ is a $k$-rational point on $X$
which is smooth on any $Y_i$ ($0 \le i \le n$).
Let us choose and fix  local parameters
 $z_1, \ldots, z_n \in
\widehat{\oo}_{x,X}$ such that
$z_{n-i+1} {|}_{Y_{i-1}} = 0$
is a local equation of the divisor
 $Y_i$
in the formal neighbourhood of the point~$x$
on the scheme $Y_{i-1}$ ($1 \le
i \le n$).  Let $\f$ be a rank~1 locally free
sheaf on $X$.
Fix a trivialization $e_x$ of the sheaf $\f$
in the  formal neghbourhood   of the point
$x$ on $X$.
Now the done choice of local parameters and trivialization
makes possible to identify
 $A_{(0,1, \ldots, n)}(\f)$
with the  $n$-dimensional local field  $k((z_1))\ldots (((z_n))$.

Moreover,
let us fix a collection of integers $0 \le j_1 \le \ldots \le j_k \le n-1$.
Define $\sigma \in S_{n-k}$ as the set
$\left\{i \: :  \: 0 \le i \le n, \, i \ne j_1, \, \ldots,
\,
i \ne j_k \right\}$.
By theorem~\ref{tend1}  we have the natural imbedding
$A_{\sigma} (\f) \lto  A_{(0,1, \ldots, n)}(\f) $.
And under identifying
of $A_{(0,1, \ldots, n)}(\f)$ with the field
$k((z_1)) \ldots ((z_n))$
the space
 $A_{\sigma} (\f)$
converts to the following
 $k$-subspace in
$k((z_1)) \ldots ((z_n))$:
\begin{equation}  \label{st}
\left\{  \sum a_{i_1, \ldots, i_n} z_1^{i_1 }
\ldots z_n^{i_n}  \; :  \;   a_{i_1, \ldots, i_n}  \in k, \,
i_{n-j_1} \ge 0,  \, i_{n-j_2} \ge 0, \, \ldots, i_{n-j_k} \ge 0
 \right\} \mbox{.}
\end{equation}

Thus, from theorem~\ref{tend2}
we obtain that for  determination  of the images of
$A_{\sigma} (\f)$
in $k((z_1)) \ldots ((z_n))$
(for any $\sigma \in S$)
it suffices to know  only one image of
$A_{(0,1, \ldots, n-1)}$
in $k((z_1)) \ldots ((z_n))$.
(All the others are obtained  by intersection
of the image of $A_{(0,1, \ldots, n-1)}$
in $k((z_1)) \ldots ((z_n))$
with the standard subspaces~(\ref{st})
in $k((z_1)) \ldots ((z_n))$.)

It is clear that these reasonings
is generalized at once to locally free sheaves
 $\f$  of rank $r$
and spaces $k((z_1)) \ldots ((z_n))^{\oplus r}$.

Denote the described map
$$
(X, Y_1, \ldots, Y_n, (z_1, \ldots, z_n), \f, e_x )
\lto A_{(0,1, \ldots, n-1)}(\f) \hookrightarrow
A_{(0,1, \ldots, n-1, n)} (\f) \stackrel{e_x}{\lto}
\qquad \qquad \qquad
$$
$$
\qquad
\qquad
\qquad
\qquad
\qquad
\qquad
\qquad
\qquad
\qquad
\quad
\stackrel{e_x}{\lto}
A_{(0,1, \ldots, n)}(\widehat{\oo}_{x,X})
\stackrel{z_1, \ldots, z_n}{\lto} k((z_1)) \ldots  ((z_n))^{\oplus r}
$$
by $\Psi_r$.

\bigskip
\bigskip
\noindent
{\bf Definition.}
$$
\begin{array}{lcl}
{\cal M}_n  & \quad \eqdef \quad &
\{X, (Y_1, \ldots, Y_n), (z_1, \ldots, z_n), \f, e_{Y_n}  \}\\[6pt]
X  && \mbox{\em
a projective equidimensional Cohen -Macaulay scheme}\\
 && \mbox{\em
of dimension $n$ over a field $k$}\\[6pt]
X=Y_0 \supset Y_1 \supset \ldots \supset Y_n
&&
\mbox{\em
a flag of closed subschemes
such that
$Y_i$ is an ample}\\ &&
\mbox{\em
 Cartier divisor on the scheme
$Y_{i-1}$ } (1 \le i \le  n)\\[6pt]
Y_n  && \mbox{\em
a smooth $k$-rational point on all $Y_i$  }  (0 \le i \le  n)\\[6pt]
z_1, \ldots, z_n &&
\mbox{\em
 formal local parameters in the point $Y_n$}\\
&& \mbox{\em such that
  $\left(z_{n-i+1} {|}_{Y_{i-1}} = 0\right) = Y_i$
in the formal}\\
&& \mbox{\em neighbourhood  of the point $Y_n$ on the
scheme $Y_{i-1}$}\\[6pt]
\f && \mbox{\em  a rank $r$ vector bundle on $X$}\\[6pt]
e_{Y_n} && \mbox{\em a trivialization of $\f$
in the formal neighbourhood}\\
&& \mbox{\em of the point $Y_n$ on $X$}
\end{array}
$$

\bigskip
In the field $K = k((z_1)) \ldots  ((z_n))$
we have the following filtration
$$
K(m) = z_n^m k((z_1)) \ldots ((z_{n-1}))[[z_n]]   \mbox{.}
$$
Let $K$-space  $V = K^{\oplus r}$, and let the filtartion
$V(m) = K(m)^{ \oplus r}$.

\bigskip
\begin{th}        \label{final}
There exists a canonical map
$$
\Phi_n \quad :
\quad
{\cal M}_n  \; \lto    \;
\left\{
\mbox{$k$-vector subspaces}
\quad
B \subset K, \quad
W \subset V
  \right\}   \mbox{}
$$
such that
\begin{enumerate}
\item   \label{f1}  from the subspace
$B \subset K$
is uniquely reconstructed
the complex $A(\oo_X)$,
which calculates cohomology of the sheaf
 $\oo_X$ on $X$;
\item  \label{f2}
from the subspace $W \subset V$
is uniquely reconstructed the complex $A(\f)$,
which calculates cohomology of the sheaf
$\f$ on $X$;
\item  \label{f3}
if $(B, W) \in \Image \Phi_n$,
then $B \cdot B \subset B$, $B \cdot W \subset W$;
\item  \label{f4}
for all $m$ the map
$$
\left\{Y_1, (Y_2, \ldots, Y_n),
(z_1, \ldots, z_{n-1}){|}_{Y_1},
\f(-mY_1){|}_{Y_1},
e_{Y_n}(-m){|}_{Y_1}\right\}
\lto
$$
$$
\lto
\left\{
\frac{B \cap K(m)}{B \cap K(m+1)}
\subset
\frac{K(m)}{K(m+1)}= k((z_1)) \ldots  ((z_{n-1}))
 \mbox{,}   \right.
$$
$$
\left.
\frac{W \cap V(m)}{W \cap V(m+1)}
\subset
\frac{V(m)}{V(m+1)}=
k((z_1)) \ldots  ((z_{n-1}))^{\oplus r}
    \right\}
$$
coincides with the map $\Phi_{n-1}$;
\item  \label{f5}
If $q, q' \in {\cal M}_n$
and $\Phi_n (q) = \Phi_n(q')$,
then $q$ is isomorphic to $q'$.
\end{enumerate}
\end{th}
\proof
If $$
q = \{ X, (Y_1, \ldots, Y_n),
(z_1, \ldots, z_n),
\f, e_{Y_n} \}  \in {\cal M}_n  \mbox{,}
$$
then to define the map $\Phi_n$
we put
$$
B = \Psi_1 (
X, Y_1, \ldots, Y_n,
(z_1, \ldots, z_n),
\oo_X,
id)
\mbox{,}  $$
$$
W = \Psi_r (
X, Y_1, \ldots, Y_n,
(z_1, \ldots, z_n),
\f,
e_{Y_n})   \mbox{,}
$$
$$
\Phi_n (q) = \left\{ B, W \right\}  \mbox{.}
$$
Now items~\ref{f1}-\ref{f4}  of this theorem
follows from theorems~\ref{teorem3}, \ref{tend1}, \ref{tend2}
and reasonings above about the map $\Psi$,
and also for item~\ref{f4} is needed the fact that
    $ (j_0)_* (j_0)^* \f = \mathop{\il}\limits_m \f(mY_1)$.

Let us show item~\ref{f5}.
Intersecting $B$
with the standard subspaces~(\ref{st}),
we can uniquely reconstruct the algebra
$A_{(0)} (\oo_X) \subset K$.
Similarly,
from $W$ we can reconstruct
the  $k$-subspace $A_{(0)} (\f) \subset V$.
Then
$$
X - Y_1 = \mathop{\rm Spec} A_{(0)} (\oo_X)
$$
\begin{equation}     \label{kvadrat}
X =
\mathop{\rm Proj} \left(
\bigoplus_{m \ge 0}
(A_{(0)} (\oo_X)  \cap K(-m))
\right)
\end{equation}
$$
\f =
\mathop{\rm Proj} \left( \bigoplus_{m \ge 0}
(A_{(0)} (\f)  \cap V(-m)) \right)   \mbox{,}
$$
where the last equalities follow from the following statement
(see~\cite[ lemma~7]{P}):
if  $X$ is a projective scheme over a field,
$\f$ is a coherent sheaf on $X$,
and $C$ is an ample divisor on $X$,
then $X \cong \mathop{\rm Proj} (S)$,
$\f \cong \mathop{\rm Proj} (F)$,
where $S = \bigoplus_{m \ge 0} H^0 (X, \oo_X (mC))$,
$F = \bigoplus_{m \ge 0} H^0 (X, \f (mC))$.

Besides, the image under the imbedding of
$$
\bigoplus_{m \ge 0}
(A_{(0)}(\oo_X) \cap K(-m+1))
\qquad
\mbox{into}   \qquad
\bigoplus_{m \ge 0}
(A_{(0)} (\oo_X)  \cap K(-m))
$$
is the homogeneous ideal
determining the subscheme
 $Y_1$ in $X$.
Now, using item~\ref{f4}  of this theorem,
we can reconstruct all the geometrical data on the subscheme
$Y_1$ in the analogous way,
and, further,  by induction
we  can reconstruct
all the data from
$q$  up to an isomorphism.
Theorem~\ref{final}
is proved.

\bigskip
\noindent
\begin{nt} {\em
Note that
$\f (nY_1) {|}_{Y_1} = \f \otimes_{\oo_X} \oo(nY_1) {|}_{Y_1}$
and the sheaf
$$
\oo(nY_1) {|}_{Y_1} = \oo(nY_1)/ \oo((n-1)Y_1) = {\cal N}_{Y_1 / X}^n
\mbox{,}$$
where the bundle ${\cal N}_{Y_1 / X}$
coincides with the normal budle of $Y_1$  in $X$
in some cases (for example, if $X$ and $Y_1$ are smooth).
}
\end{nt}

\begin{nt} {\em
From~(\ref{kvadrat})
and absence of  divisors of zero in the field $K$
it follows at once that the schemes $X, Y_1, \ldots, Y_n$,
satisfied the conditions of the definiton ${\cal M}_n$,
are always integral schemes.
}
\end{nt}

\begin{nt} {\em
 $\Phi_1$
is a variant of the Krichever correspondence for curves
 (see~\cite{M}, \cite{P}).
Besides,
any integral noetherian scheme of dimension $1$
is a Cohen-Macaulay scheme.

$\Phi_2$
coincides with the map, constructed in~\cite{P}.
Note that any normal noetherian scheme of dimension
$2$
is a Cohen-Macaulay  scheme (see~\cite[ ch.~II, th.~8.22A]{Ha}).
Also in~\cite{P}
is analyzed an example of the Krichever map for $X = {\bf P}_2$.
}
\end{nt}

\begin{nt}
{\em
In~\cite{P}
is discussed the problem of change
of locally free sheaves to torsion free sheaves.

Let $X$ be a smooth projective surface
with a flag of irreducible subvarieties
 $Y_1 \supset Y_2$
such that $Y_1$ is an ample divisor on $X$, and $Y_2$ is a point.
Let $\g$ be any torsion free sheaf on $X$.
Then we have the imbedding
\begin{equation}     \label{posledo}
0 \lto \g \lto \g^{**}   \mbox{.}
\end{equation}
The sheaf $\g^{**}$
is reflexive; and since $\dm X =2$,
we have that $\g^{**}$ is a locally free sheaf.
Now applying the exact functors $A_{\sigma}$
and the left exact functor  $H^0 (X, \cdot)$
to sequence~(\ref{posledo}),
we obtain at once from the obtained sequences
and theorem~\ref{tend1}
for the sheaf $\g^{**}$ that theorem~\ref{tend1}
is true for the sheaf $\g$ as well.

But theorem~\ref{tend2}  is no longer true
for torsion free sheaves on $X$.
Indeed, let $\g = m_Q$,
where $m_Q$ is the ideal sheaf of a point $Q \in X$.
Then, applying the exact functors $A_{\sigma}$
to the exact sequence
$$
0 \lto m_Q \lto \oo_X \lto k_Q \lto 0
$$
we obtain at once that
\begin{itemize}
\item if $Q \notin Y_1$,
then
$$
A_{(0)} (\g) \ne A_{(01)} (\g)  \cap A_{(02)} (\g) \mbox{,}
$$
because $A_{(01)} (\g)  \ne A_{(01)} (\oo_X)$,
but $A_{(01)} (\g) = A_{(01)} (\oo_X)$,
$A_{(02)} (\g) = A_{(02)} (\oo_X)$;
\item
if $Q \in Y_1$, $Q \ne Y_2$,
then
$$
A_{(1)}(\g) \ne A_{(01)} (\g) \cap A_{(12)} (\g) \mbox{,}
$$
because $A_{(1)} (\g) \ne A_{(1)} (\oo_X)$,
but $A_{(01)} (\g) = A_{(01)} (\oo_X)$,
$A_{(12)} (\g) = A_{(12)} (\oo_X)$.
\end{itemize}
}
\end{nt}

\bigskip
%******************************************************************
%******************************************************************

\nopagebreak[1]
\bigskip
\bigskip
 Steklov Mathematical Institute
\medskip \\
 {\em E-mail :} d\_osipov@mi.ras.ru \\
 \hphantom{{\em E-mail :}} d\_osipov@mail.ru

%********************************************************************
\end{document}